\PassOptionsToPackage{unicode}{hyperref}
\PassOptionsToPackage{hyphens}{url}
\documentclass[
]{article}
\usepackage{amsmath,amssymb}
\usepackage{iftex}
\ifPDFTeX
  \usepackage[T1]{fontenc}
  \usepackage[utf8]{inputenc}
  \usepackage{textcomp} 
\else 
  \usepackage{unicode-math} 
  \defaultfontfeatures{Scale=MatchLowercase}
  \defaultfontfeatures[\rmfamily]{Ligatures=TeX,Scale=1}
\fi
\usepackage{lmodern}
\ifPDFTeX\else
\fi
\IfFileExists{upquote.sty}{\usepackage{upquote}}{}
\IfFileExists{microtype.sty}{
  \usepackage[]{microtype}
  \UseMicrotypeSet[protrusion]{basicmath} 
}{}
\makeatletter
\@ifundefined{KOMAClassName}{
  \IfFileExists{parskip.sty}{%
    \usepackage{parskip}
  }{
    \setlength{\parindent}{0pt}
    \setlength{\parskip}{6pt plus 2pt minus 1pt}}
}{
  \KOMAoptions{parskip=half}}
\makeatother
\usepackage{xcolor}
\usepackage{longtable,booktabs,array}
\usepackage{calc} 
\usepackage{etoolbox}
\makeatletter
\patchcmd\longtable{\par}{\if@noskipsec\mbox{}\fi\par}{}{}
\makeatother
\IfFileExists{footnotehyper.sty}{\usepackage{footnotehyper}}{\usepackage{footnote}}
\makesavenoteenv{longtable}
\usepackage{graphicx}
\makeatletter
\def\maxwidth{\ifdim\Gin@nat@width>\linewidth\linewidth\else\Gin@nat@width\fi}
\def\maxheight{\ifdim\Gin@nat@height>\textheight\textheight\else\Gin@nat@height\fi}
\makeatother
\setkeys{Gin}{width=\maxwidth,height=\maxheight,keepaspectratio}
\makeatletter
\def\fps@figure{htbp}
\makeatother
\providecommand{\tightlist}{%
  \setlength{\itemsep}{0pt}\setlength{\parskip}{0pt}}
\usepackage{amsthm}
\newtheoremstyle{paperplain}
  {6pt}
  {6pt}
  {\itshape}
  {}
  {}
  {}
  {0pt}
  {}
\theoremstyle{paperplain}
\newtheorem*{papertheorem}{}
\newtheorem*{paperlemma}{}
\newtheorem*{paperproposition}{}
\newtheorem*{papercorollary}{}

\newtheorem*{paperclaim}{}

\ifLuaTeX
  \usepackage{selnolig}  
\fi
\usepackage[numbers,square]{natbib}
\IfFileExists{bookmark.sty}{\usepackage{bookmark}}{\usepackage{hyperref}}
\IfFileExists{xurl.sty}{\usepackage{xurl}}{} 
\hypersetup{
  pdftitle={Interchange graphs of (0,1)-matrices are maximally Hamiltonian},
  pdfauthor={Jeffrey S. Baggett; Huiya Yan},
  hidelinks,
  pdfcreator={LaTeX via pandoc}}

\title{Interchange graphs of \((0,1)\)-matrices are maximally
Hamiltonian}
\author{Jeffrey S. Baggett\footnote{Department of Mathematics and
  Statistics, University of Wisconsin--La Crosse. Email:
  jbaggett@uwlax.edu. Corresponding author.} \and Huiya Yan\footnote{Department
  of Mathematics and Statistics, University of Wisconsin--La Crosse.
  Email: hyan@uwlax.edu.}}
\date{}

\begin{document}
\maketitle

\begin{abstract}

For integer vectors \(R,S\) let \(\mathcal A(R,S)\) denote the class of
\((0,1)\)-matrices with row sum vector \(R\) and column sum vector
\(S\). Its \emph{interchange graph} \(G(R,S)\) has \(\mathcal A(R,S)\)
as its vertex set, two matrices being adjacent when they differ by a
single \(2\times 2\) interchange. Brualdi asked whether \(G(R,S)\) is
Hamiltonian for every \(R,S\). We prove the stronger statement that
\(G(R,S)\) is \textbf{maximally Hamiltonian}: Hamilton-laceable when
bipartite, and Hamilton-connected when not. The proof is a structural
induction on the number of matrices in the class, organized by the
structure theory of interchange graphs. Deleting inactive lines and
splitting invariant positions expresses any class as a Cartesian
product, reducing the argument to the prime factors. The bipartite
classes are products of complete transposition graphs; we settle them
together, without induction, by proving they are paired
\(2\)-disjoint-path-coverable and hence Hamilton-laceable, using a
recent theorem of Coleman, Fischberg, Gong, Harrington and Wong on
paired disjoint path covers. The non-bipartite classes divide into three
cases: products assembled from smaller factors, a base of Johnson graphs
and small classes, and the large prime classes, treated by a
pivot-and-fiber construction whose line quotients are matroid
base-exchange graphs. The complete argument has been machine-checked in
the Lean 4 proof assistant from first principles together with seven
cited results of the literature; the disjoint-path-cover results it
imports are themselves proved within the formalization.

\end{abstract}

(secondary). \emph{Keywords:} interchange graph, (0,1)-matrix,
Hamilton-connected, Hamilton-laceable, disjoint path cover,
transportation polytope, combinatorial Gray code.

\begin{center}\rule{0.5\linewidth}{0.5pt}\end{center}

\hypertarget{introduction}{%
\subsection{1. Introduction}\label{introduction}}

Let \(R=(r_1,\dots,r_m)\) and \(S=(s_1,\dots,s_n)\) be nonnegative
integer vectors with \(\sum r_i=\sum s_j\), and let \(\mathcal A(R,S)\)
be the set of \(m\times n\) \((0,1)\)-matrices with row sums \(R\) and
column sums \(S\). A \textbf{\(2\times 2\) interchange} replaces a
submatrix \(\left[\begin{smallmatrix}1&0\\0&1\end{smallmatrix}\right]\)
by \(\left[\begin{smallmatrix}0&1\\1&0\end{smallmatrix}\right]\), or
vice versa, leaving all other entries fixed. Ryser's classical theorem
\citep{ryserCombinatorialPropertiesMatrices1957} states that any two
members of \(\mathcal A(R,S)\) are connected by a sequence of
interchanges. The \textbf{interchange graph} \(G(R,S)\) records this:
its vertices are the matrices of \(\mathcal A(R,S)\), and two are
adjacent precisely when one is obtained from the other by a single
interchange. Interchange graphs are central objects of combinatorial
matrix theory; they are spanning subgraphs of the \(1\)-skeleta of
transportation polytopes (every interchange is a polytope edge, though
the skeleton has further edges, along longer alternating cycles), and
they are the state graphs of the standard swap Markov chains on
contingency tables.

Brualdi \citep[Problem 3.7]{brualdiMatricesZerosOnes1980} posed the
following problem, stated there in the matrix language recalled above.

\begin{quote}
\textbf{Problem (Brualdi, 1980, Problem 3.7).} Does the interchange
graph \(G(R,S)\) have a Hamilton cycle? If it has no Hamilton cycle,
does it have a Hamilton path? If it has no Hamilton path, is there a
matrix \(A\in\mathcal A(R,S)\) and a family of paths from \(A\) covering
every other matrix exactly once?
\end{quote}

The question leads with Hamiltonicity, asking the weaker path and
path-cover properties only as successive fallbacks. It is one instance
of a general principle, phrased by Mütze
\citep{mutzeCombinatorialGrayCodes2023} as a meta-conjecture, that a
flip graph is Hamiltonian once the evident connectivity, minimum-degree,
and balancedness obstructions are absent. The problem is known only in
special cases. Hamilton cycles were obtained for certain special margin
vectors by Zhang--Zhang and H. Zhang (see Brualdi
\citep[§6.3]{brualdiCombinatorialMatrixClasses2006} for an account). Li
and Zhang \citep[Thm. 2.3,
pp.~109--110]{liHamiltonicityTypeInterchange1994a} proved Hamiltonicity
when one margin is all-ones (equivalently, by complementation,
all-\((\max-1)\)). In the paragraph following that proof they state that
the same method gives edge-Hamiltonicity \citep[
p.~110]{liHamiltonicityTypeInterchange1994a}. This class consists of
labeled set partitions with prescribed part sizes and is more general
than the permutation-margin case; the complete transposition graph
occurs only when both margins are all-ones. A \textbf{split graph} is a
graph whose vertex set partitions into a clique and an independent set,
the edges between the two parts being unrestricted. Arikati and Peled
\citep[pp.~214--215]{arikatiRealizationGraphDegree1999} gave the exact
split-graph translation of a bipartite realization graph described
below, and established Hamiltonicity for degree sequences of
majorization gap one. This realization-graph viewpoint is surveyed by
Barrus \citep{barrusRealizationGraphsDegree2016}, where the general
Hamiltonicity question is likewise recorded as open. The general problem
has remained open; it appears as Problem P60 in Mütze's survey
\citep{mutzeCombinatorialGrayCodes2023}. Independently, and at the same
time as the present work, Hladík and Fink
\citep{hladikRealizationGraphEvery2026} announce, in a preprint, that
the realization graph of every degree sequence has a Hamilton path
starting from any prescribed vertex; the interchange graphs are the
bipartite case of that setting. We prove the stronger
maximal-Hamiltonicity property defined below; in particular, every
interchange graph with at least three matrices has a Hamilton cycle, the
first, Hamilton-cycle form of Brualdi's problem.

We establish a strengthening. Recall that a connected graph is
\textbf{Hamilton-connected} if every pair of distinct vertices is joined
by a Hamilton path, and that a connected bipartite graph is
\textbf{Hamilton-laceable} if every pair of vertices in different color
classes is joined by a Hamilton path (for a bipartite graph the latter
is the strongest possible such property). Call a graph \textbf{maximally
Hamiltonian} if it is Hamilton-laceable when bipartite and
Hamilton-connected when not, a dichotomy with precedent in Chen and
Quimpo's work on abelian group graphs
\citep{chenStronglyHamiltonianAbelian1981}.

\begin{papertheorem}

\textbf{Theorem 1.1.} For all realizable \(R,S\) (that is, with
\(\mathcal A(R,S)\ne\emptyset\)), the interchange graph \(G(R,S)\) is
maximally Hamiltonian. Consequently, if \(|\mathcal A(R,S)|\ge3\), then
\(G(R,S)\) has a Hamilton cycle; if \(|\mathcal A(R,S)|=2\), the two
matrices give Brualdi's cyclic Gray listing.

\end{papertheorem}

Maximal Hamiltonicity, rather than bare Hamiltonicity, is the property
our induction can carry: a Hamilton cycle is not preserved under the
Cartesian-product and fiber gluings that organize the argument, whereas
Hamilton-connectedness and -laceability are. For a graph of order at
least three, take any edge, join its endpoints by the Hamilton path
supplied by maximal Hamiltonicity, and close with that edge. Thus every
edge lies on a Hamilton cycle. The two-vertex interchange graph is
\(K_2\); its sole edge gives the two-term cyclic listing, although
\(K_2\) is not a cycle under the standard simple-graph convention. The
one-vertex class is maximally Hamiltonian only vacuously and has no
Hamilton cycle under that convention. Accordingly, a statement of
Brualdi's problem that includes the singleton must either exclude it
from the Hamilton-cycle conclusion or declare its one-term listing to be
a degenerate cycle. The Lean companion also covers the empty class
vacuously, so its statement is formally stronger than the display above.

\begin{papercorollary}

\textbf{Corollary.} The realization graph of any realizable bipartite
degree sequence (all labeled bipartite graphs with prescribed degrees on
each side, adjacent under 2-switches) is maximally Hamiltonian.
Equivalently, by the correspondence of Arikati--Peled, so is the
realization graph of any split-graph degree sequence --- the degree
sequence of a graph that partitions into a clique and an independent
set. (The analogous question for arbitrary degree sequences, cf.~Barrus,
remains open.)

\end{papercorollary}

\emph{Proof.} A \((0,1)\)-matrix with margins \(R,S\) is the biadjacency
matrix of a labeled bipartite graph on fixed sides
\(X=\{x_1,\ldots,x_m\}\) and \(Y=\{y_1,\ldots,y_n\}\), and an
interchange is exactly a color-class-preserving \(2\)-switch. This
identifies its bipartite realization graph with \(G(R,S)\).

For the split correspondence, add every edge inside \(X\) and no edge
inside \(Y\). A bipartite realization with side degrees \(R=(r_i)\) and
\(S=(s_j)\) becomes a split realization on the labeled partition
\(X\mathbin{\dot\cup}Y\) with degree sequence \[
d_{R,S}=(m-1+r_1,\ldots,m-1+r_m,s_1,\ldots,s_n).
\] Deleting the fixed clique edges reverses the map. A \(2\)-switch in a
split realization preserves its clique--independent-set partition. For a
\(2\)-switch is admissible only when the two edges it creates are
currently absent, and since \(Y\) is independent every edge of a split
realization has an end in \(X\); so deleting a clique edge would, on
either pairing with its partner, create a second edge with both ends in
\(X\), and every such pair is already present. Hence no \(2\)-switch
touches a clique edge, and the same observation leaves, for a pair of
crossing edges, only the pairing that produces two crossing edges --- a
bipartite \(2\)-switch after the clique edges are deleted. Connectivity
then supplies surjectivity: any two realizations of a degree sequence
are joined by \(2\)-switches, so every labeled realization of
\(d_{R,S}\) is reachable from the displayed one, and by the previous
sentence every step retains the same labeled partition. This is the
bijection onto all labeled realizations of \(d_{R,S}\) observed by
Arikati and Peled, hence an isomorphism of realization graphs
\citep[pp.~214--215]{arikatiRealizationGraphDegree1999}. Conversely,
choose the labeled clique \(X\) and independent set \(Y\) of a split
sequence \(d\); subtract \(m-1\) from the degrees on \(X\) and retain
the degrees on \(Y\) to recover \((R,S)\) and the inverse isomorphism.
Theorem 1.1 now applies in both directions. \(\square\)

\textbf{Outline.} Section 2 fixes notation and recalls the structure
theory of interchange graphs. The proof of Theorem 1.1 is a structural
induction on the number of matrices in the class; Section 3 lays out its
architecture. The primary split is on whether \(G(R,S)\) is bipartite.

A bipartite class is a Cartesian product of complete transposition
graphs (Proposition 2.2), and we settle every such class at once,
without induction, through the following paired-cover property.

\begin{paperlemma}

\textbf{Lemma 1.2 (Bipartite 2-DPC Lemma).} Every balanced bipartite
interchange graph is paired \(2\)-disjoint-path-coverable.

\end{paperlemma}

Here a \emph{balanced bipartite} interchange graph is one that is
bipartite with color classes of equal size; the
paired-disjoint-path-cover property is defined in Section 2.4. A paired
cover yields a Hamilton path between any two opposite-colored vertices,
so the Bipartite 2-DPC Lemma gives Hamilton-laceability for the whole
bipartite branch. Section 7 proves it: by a classical theorem of Brualdi
every balanced bipartite interchange graph is a Cartesian product of
complete transposition graphs, and an induction on the number of
factors, built on the Coleman--Fischberg--Gong--Harrington--Wong
theorem, shows every such product is paired
\(2\)-disjoint-path-coverable.

A non-bipartite class is treated by the induction and shown
Hamilton-connected. After deleting inactive lines --- rows and columns
of sum \(0\) or maximal, whose entries are therefore the same in every
matrix of the class, so that no interchange meets them --- a surviving
invariant position (Section 2.2), meaning a row--column position whose
entry likewise never varies, expresses the class as a nontrivial
Cartesian product of smaller interchange graphs, and Section 4 lifts
maximal Hamiltonicity from the factors. A class with no invariant
position is prime: its base cases, the Johnson graphs \(J(f,k)\) (the
graph on the \(k\)-subsets of an \(f\)-set, two subsets adjacent when
they share \(k-1\) elements) and the classes of order at most six, are
handled in Section 6, and the remaining large prime classes by the
pivot-and-fiber construction of Section 5: fix a single row, group the
matrices of the class by what that row does --- the groups are its
\emph{fibers} --- and build the path fiber by fiber. The quotient that
records how the fibers meet is a matroid base-exchange graph. The
Bipartite 2-DPC Lemma returns as a subroutine inside Section 4 whenever
a non-bipartite product has a bipartite factor.

\textbf{Verification.} The argument has been checked twice by machine,
in complementary ways. The entire proof --- the reduction of Sections
3--6, the pivot construction of Section 5, and the factor induction of
Section 7 --- has been formalized in the Lean 4 proof assistant, so the
composed statement, that every interchange graph is maximally
Hamiltonian, is machine-checked from Lean's foundations together with
exactly seven cited results of the refereed literature (six named
theorems and a classical transportation-polytope dimension count, listed
with their sources and roles in Table 1). The disjoint-path-cover
results of Section 7 that the argument imports from Coleman et al.
\citep{colemanPairedn1ton1Disjoint2025} and Jo, Park and Chwa
\citep{joPaired2disjointPath2013} are themselves proved within the
formalization from first principles, so those citations record
provenance rather than assumption. Independently, the theorem was
confirmed on all small classes by an exhaustive SAT-based computation.
Neither check replaces the cited results, which are used as published,
and neither forms part of the argument: the human-readable proof of
Sections 3--7 stands on its own, and the mechanization and the SAT sweep
corroborate it. Section 9 records both --- the release artifact, the
build gate and axiom trace, and the claim-by-claim correspondence,
together with the computational sweep.

\begin{center}\rule{0.5\linewidth}{0.5pt}\end{center}

\hypertarget{preliminaries}{%
\subsection{2. Preliminaries}\label{preliminaries}}

Throughout, \(G(R,S)\) denotes the interchange graph of
\(\mathcal A(R,S)\) --- which we call the \textbf{class} --- \(\square\)
the Cartesian product of graphs, and \(\chi\) the interchange-parity
\(2\)-coloring of a bipartite class (Section 2.3). The two
normalizations that recur, deleting inactive lines and splitting
invariant positions, are collected in Section 2.2.

\hypertarget{lines-fibers-and-quotients}{%
\subsubsection{2.1 Lines, fibers, and
quotients}\label{lines-fibers-and-quotients}}

A \textbf{line} of a matrix is a row or a column; its \textbf{pattern}
is its support, the set of positions in which it carries a \(1\). Fix a
line \(L\). For each realizable pattern \(p\) of \(L\), the
\textbf{fiber} \(F_p\) is the set of matrices of \(\mathcal A(R,S)\)
whose line \(L\) equals \(p\); the fibers partition \(\mathcal A(R,S)\)
according to the value of \(L\). All matrices in a fiber \emph{share}
the same \(L\), so a single \(2\times2\) interchange joins two of them
only when it avoids \(L\) entirely (an interchange meeting \(L\) would
change its pattern, moving to a different fiber). Deleting the common
line \(L\) therefore identifies \(F_p\) with the interchange graph of
the smaller class on the remaining lines (the residual class obtained by
removing \(L\) and reducing the margins by \(p\)). The interchanges that
\emph{do} meet \(L\) pass between fibers and are recorded by the
\textbf{quotient} \(Q_L\): its vertices are the realizable patterns of
\(L\), with \(p,q\) adjacent when a single interchange (necessarily
meeting \(L\)) carries a matrix of \(F_p\) to one of \(F_q\);
equivalently, when \(q\) is obtained from \(p\) by relocating one entry
of \(L\). All of these notions are illustrated and used in Section 5.

\hypertarget{invariant-positions}{%
\subsubsection{2.2 Invariant positions}\label{invariant-positions}}

A \textbf{cell} is a row--column position of the matrices of the class,
that is, a single matrix entry \((i,j)\). A cell is \textbf{invariant}
if it takes the same value in every matrix of \(\mathcal A(R,S)\); such
a cell is incident to no interchange.

Invariant positions arise in two ways, which must be separated. An
\textbf{inactive line} (a row or column whose sum is \(0\) or maximal)
consists entirely of constant cells; iteratively deleting inactive lines
until none remain leaves a class with the \emph{same} interchange graph
in which every line is \textbf{active} (\(0<r_i<n\) and \(0<s_j<m\)). We
assume throughout that the class is active; this deletion is
graph-isomorphic, not a factorization. For a class with active lines,
the \emph{remaining} invariant positions induce a genuine Cartesian
factorization: after permuting rows and columns, the non-invariant cells
occupy disjoint diagonal blocks, every off-diagonal rectangle forced to
all-ones or all-zeros, and the matrices vary independently within the
blocks. Brualdi's display identifies the class with the product of the
block classes and the interchange graph with the Cartesian product of
their interchange graphs \citep[§6.1,
pp.~285--286]{brualdiCombinatorialMatrixClasses2006}. The same
construction appears in Brualdi--Manber: their Theorem 1 gives the
invariant block form, and the paragraph following it identifies the two
graph factors \citep[Thm. 1 and following paragraph,
pp.~158--159]{brualdiPrimeInterchangeGraphs1983}. Thus every factor is
itself the interchange graph of a matrix block, not merely an abstract
Sabidussi factor. In an active class both displayed factors have at
least two vertices
\citep[pp.~158--159]{brualdiPrimeInterchangeGraphs1983}; consequently
each has order strictly smaller than their product. Brualdi--Manber's
Theorem 9 gives the converse and the exact prime criterion \citep[Thm.
9, p.~169]{brualdiPrimeInterchangeGraphs1983}: an active class is
\textbf{prime} (admits no nontrivial Cartesian factorization) exactly
when it has \textbf{no invariant positions}. We call such a class
\textbf{invariant-free}. (Primeness here means
Cartesian-indecomposability, not connectivity: \(C_4=K_2\,\square\,K_2\)
is connected yet decomposable.) Thus, once inactive lines are deleted,
every interchange graph is a Cartesian product of invariant-free (prime)
factors: the reduction of Section 3 splits this product off first
(Section 4) and otherwise treats an invariant-free class. (Fixing a line
pattern to form a fiber, as in Section 5, can create new inactive lines
or invariant positions; the normalization (delete inactive lines, then
split off the product) is re-applied at each recursive step, the
invariant carried by the induction being ``strictly smaller interchange
graph''.)

\hypertarget{bipartiteness}{%
\subsubsection{2.3 Bipartiteness}\label{bipartiteness}}

\begin{paperlemma}

\textbf{Lemma 2.1} (Brualdi \citep[Thm. 6.3.4(i)\(\Leftrightarrow\)(ii),
p.~298]{brualdiCombinatorialMatrixClasses2006}). \emph{An invariant-free
interchange graph is bipartite if and only if it is triangle-free.
Consequently, an arbitrary interchange graph is bipartite if and only if
it is triangle-free.}

\end{paperlemma}

The cited theorem gives the first sentence, which is the direction used
in the prime branch. For the consequence, delete inactive lines and use
the block product of Section 2.2. A Cartesian product is bipartite
exactly when every nontrivial factor is bipartite, and it contains a
triangle exactly when one factor contains a triangle. Applying the
invariant-free statement factor by factor proves the general form. This
equivalence is special to interchange graphs; triangle-free does not
imply bipartite for a general graph.

When \(G(R,S)\) is bipartite, fix any \(A_0\in\mathcal A(R,S)\) and
assign to each matrix the parity of the length of an interchange path
from \(A_0\) to it. All such paths have the same parity because the
graph is bipartite, so this is well defined, and the two resulting
classes are the color classes. We call it the \textbf{interchange
parity} and write \(\chi\) for the corresponding \(2\)-coloring.

\hypertarget{disjoint-path-covers}{%
\subsubsection{2.4 Disjoint path covers}\label{disjoint-path-covers}}

Let \(H\) be a graph and let \(k\ge1\). A \textbf{\(k\)-demand} is a
family of \(k\) terminal pairs \((a_1,b_1),\ldots,(a_k,b_k)\) whose
\(2k\) terminals are distinct. A \textbf{paired \(k\)-disjoint path
cover} for such a demand is a family of \(k\) pairwise vertex-disjoint
paths, the \(i\)th joining \(a_i\) to \(b_i\), whose vertex sets
together exhaust \(V(H)\). We abbreviate a \(2\)-demand to
\textbf{demand} and a paired \(k\)-disjoint path cover to a
\textbf{paired \(k\)-cover}, since \(k=2\) is the case the argument
uses; the general \(k\) is needed only to state the two theorems of
Section 7.1, which produce covers of higher order and then descend to
\(k=2\).

For a bipartite graph with color classes \(V_1,V_2\), a \(k\)-demand is
\textbf{pairwise-opposite} if each of its pairs joins the two classes;
equivalently, after relabeling within pairs (paths are undirected),
\(a_i\in V_1\) and \(b_i\in V_2\) for every \(i\). For \(k=2\) this
reads \(\{a_1,a_2\}\subseteq V_1\) and \(\{b_1,b_2\}\subseteq V_2\). A
bipartite graph is \textbf{paired \(2\)-disjoint-path-coverable} (paired
\(2\)-DPC) if it admits a paired \(2\)-disjoint path cover for
\emph{every} pairwise-opposite demand. This property is strictly
stronger than Hamilton-laceability, its one-pair case.

Parts of the disjoint-path-cover literature
\citep{colemanPairedn1ton1Disjoint2025, joPaired2disjointPath2013}
instead use the wider \textbf{balanced-union} convention, which for two
pairs also allows one same-color pair in \(V_1\) and one same-color pair
in \(V_2\). Every theorem we invoke has the special form
\(S\subseteq V_1\), \(T\subseteq V_2\), so every demand arising here is
pairwise-opposite, and we use that term throughout rather than
``balanced'' --- which in this paper is reserved for a bipartite graph
whose two color classes have equal size.

\hypertarget{complete-transposition-graphs}{%
\subsubsection{2.5 Complete transposition
graphs}\label{complete-transposition-graphs}}

The \textbf{complete transposition graph} \(CT_a\) is the Cayley graph
of the symmetric group \(S_a\) with respect to the set of all
transpositions. It is exactly the interchange graph of the all-ones
margins \(R=S=(1,\dots,1)\) on \(a\) lines; thus \(|V(CT_a)|=a!\). It is
bipartite (transpositions are odd permutations), vertex-transitive, and
(for \(a\ge2\)) balanced. The complementary margins
\(R=S=(a-1,\dots,a-1)\) give the \emph{co-permutation} class, which is
isomorphic to \(CT_a\) via the complementation \(A\mapsto J-A\), where
\(J\) is the \(a\times a\) all-ones matrix (an isomorphism of
interchange graphs).

\begin{paperproposition}

\textbf{Proposition 2.2} (synthesis: Brualdi \citep[§6.1, pp.~285--286,
§6.3, p.~295, Thm. 6.3.4,
p.~298]{brualdiCombinatorialMatrixClasses2006}, Brualdi--Manber
\citep[Thms. 1 and 9, pp.~158--159,
169]{brualdiPrimeInterchangeGraphs1983}, and Sabidussi unique
factorization \citep{sabidussiGraphMultiplication1959}; see proof).
\emph{Every bipartite interchange graph is a Cartesian product of
complete transposition graphs and their co-permutation copies; in
particular, when nontrivial (order \(\ge 2\)), it is \textbf{balanced}
and of even order.} (Balance is thus automatic, not an extra hypothesis:
each prime bipartite block is a complete-transposition or co-permutation
graph, which is balanced, so any product of them is.)

\end{paperproposition}

\emph{Proof.} For an invariant-free class, Theorem
6.3.4(ii)\(\Leftrightarrow\)(iii) of
\citep[p.~298]{brualdiCombinatorialMatrixClasses2006} states that
triangle-freeness (equivalently, by Lemma 2.1, bipartiteness) holds
precisely when \(m=n\) and either \(R=S=(1,\dots,1)\) (a permutation
block, \(CT_n\)) or \(R=S=(n-1,\dots,n-1)\) (a co-permutation block). A
decomposable bipartite class factors through the invariant blocks of
Section 2.2 \citetext{\citealp[§6.1,
pp.~285--286]{brualdiCombinatorialMatrixClasses2006}; \citealp[Thm. 1
and following paragraph,
pp.~158--159]{brualdiPrimeInterchangeGraphs1983}}. Sabidussi's unique
Cartesian factorization \citep{sabidussiGraphMultiplication1959}
identifies the resulting indecomposable factors, and Brualdi--Manber's
Theorem 9 \citep[p.~169]{brualdiPrimeInterchangeGraphs1983} says that
each such active matrix block is invariant-free. Each factor is also
triangle-free, hence by the above is a permutation or co-permutation
block. The general bipartite case reduces to the invariant-free core
(§2.2). \(\square\)

\begin{center}\rule{0.5\linewidth}{0.5pt}\end{center}

\hypertarget{the-reduction}{%
\subsection{3. The reduction}\label{the-reduction}}

Theorem 1.1 is proved by induction on \(|V(G)|\), with the Bipartite
2-DPC Lemma (Lemma 1.2) as a global hypothesis (never carried
inductively; discharged independently in Section 7). The induction
hypothesis is that every strictly smaller interchange graph is maximally
Hamiltonian. We split into cases, each reducing \(G\) either to strictly
smaller interchange graphs or to the Bipartite 2-DPC Lemma. The
following map records the branches; the paragraphs below carry them out.

\textbf{Induction bases.} If \(|V(G)|=0\), maximal Hamiltonicity holds
vacuously; this case occurs only in the stronger formal statement, since
Theorem 1.1 assumes realizability. If \(|V(G)|=1\), it again holds
vacuously. These two cases are removed before we invoke balance,
pairwise-opposite demands, or Cartesian prime factorization. Hence
assume \(|V(G)|\ge2\) for the remainder of the reduction.

\begin{quote}
\textbf{Proof architecture.} Normalize \(G(R,S)\) by deleting inactive
lines. If \(G\) is \emph{bipartite} it is a Cartesian product of
complete transposition graphs (Proposition 2.2), paired
\(2\)-disjoint-path-coverable by the Bipartite 2-DPC Lemma (Section 7)
and hence Hamilton-laceable, with no appeal to induction. Otherwise
\(G\) is \emph{non-bipartite} and Hamilton-connected by induction on
\(|V(G)|\), in three cases: (i) an invariant position remains, so
\(G=A\,\square\,B\) is a nontrivial product of smaller interchange
graphs, lifted in Section 4; (ii) \(G\) is prime and a base case, a
Johnson graph or a class of order at most \(6\), treated in Section 6;
(iii) \(G\) is prime with at least three active rows and columns and
\(|V(G)|>6\), treated by the pivot-and-fiber construction of Section 5,
whose line quotients are matroid base-exchange graphs.
\end{quote}

\textbf{Bipartite classes.} If \(G\) is bipartite (Lemma 2.1), then by
Proposition 2.2 it is a Cartesian product of complete-transposition and
co-permutation blocks; by the Bipartite 2-DPC Lemma (Section 7) it is
paired \(2\)-disjoint-path-coverable, hence Hamilton-laceable by Theorem
7.3 with \(\ell=1\) whenever \(|V(G)|\ge4\)
\citep{colemanPairedn1ton1Disjoint2025}; the two-vertex class \(K_2\) is
Hamilton-laceable directly. This settles every bipartite class at once
(single blocks included), and uses no induction.

\textbf{Non-bipartite classes.} First delete inactive lines (§2.2) to
reach the graph-isomorphic active core; assume \(G\) is so normalized,
and non-bipartite.

\begin{itemize}
\tightlist
\item
  If the active class still has an invariant position, then by §2.2
  \(G\) is a nontrivial Cartesian product \(A\,\square\,B\) of strictly
  smaller active interchange graphs, at least one of which is
  non-bipartite (a product is bipartite iff both factors are). Section 4
  lifts maximal Hamiltonicity from the factors (induction hypothesis) to
  \(G\).
\item
  Otherwise the active class has no invariant position, hence is
  \textbf{prime} (invariant-free; Brualdi--Manber, §2.2). If at most two
  rows are active, or at most two columns are active, \(G\) is a Johnson
  graph \(J(f,k)\), handled in Section 6 (Alspach). Otherwise at least
  three rows and at least three columns are active; if \(|V(G)|\le 6\)
  it is a base class (Section 6), and if \(|V(G)|>6\) it is
  indecomposable non-base of the genuinely many-line kind treated in
  Section 5.
\end{itemize}

These cases are exhaustive, and each yields maximal Hamiltonicity of
\(G\): Hamilton-laceability when bipartite, Hamilton-connectedness
otherwise. Sections 4--6 carry out the non-bipartite cases; Section 7
proves the Bipartite 2-DPC Lemma.

\textbf{Well-foundedness and dependencies.} The induction runs on
\(|V(G)|\), and every non-bipartite case reduces \(G\) to strictly
smaller interchange graphs or to the Bipartite 2-DPC Lemma, so no case
appeals to \(G\) itself. The Bipartite 2-DPC Lemma is not part of this
induction. Section 7 proves it separately, by an induction on the number
of Cartesian factors resting on the disjoint-path-cover theorem of
Coleman, Fischberg, Gong, Harrington and Wong, and the bipartite case
here and the one-bipartite-factor case of Section 4 invoke it as a
completed fact. Within the main induction, Section 4 applies the
hypothesis to the factors \(A,B\); Section 5 applies the induction
hypothesis to the fibers of a \emph{separating} line --- one on which
the two prescribed endpoints have different patterns (§5) --- and, after
proving its matroid base-exchange quotient non-bipartite, invokes the
Naddef--Pulleyblank dichotomy to obtain Hamilton-connectedness; Section
6 is the base, proved directly. The dependencies form no cycle: Section
7 stands alone, Sections 4 and 5 descend in \(|V(G)|\), and Section 6
terminates the recursion.

\begin{center}\rule{0.5\linewidth}{0.5pt}\end{center}

\hypertarget{cartesian-products}{%
\subsection{4. Cartesian products}\label{cartesian-products}}

This section handles case (i) of the non-bipartite branch, the
decomposable classes. It shows that maximal Hamiltonicity lifts across a
Cartesian product, given it for the two factors by the induction
hypothesis. The one nontrivial input is the Bipartite 2-DPC Lemma of
Section 7, needed when a factor is bipartite; the output,
Hamilton-connectedness of the product, returns to the induction of
Section 3.

Let \(G=A\,\square\,B\) be the nontrivial Cartesian product produced by
the invariant-position decomposition (§2.2). Both \(A\) and \(B\) are
strictly smaller interchange graphs, so the strong induction hypothesis
makes each of them maximally Hamiltonian. As \(G\) is non-bipartite, at
least one factor is non-bipartite; we show \(G\) is Hamilton-connected.

View \(G\) as copies (``layers'') of \(A\) indexed by \(V(B)\): layers
over adjacent \(B\)-vertices are joined by the identity matching on the
\(A\)-coordinate. Fix distinct endpoints \(s_0,s_1\). We build a
Hamilton \(s_0\)--\(s_1\) path by following a spanning walk \(W\) of
\(B\) from the layer of \(s_0\) to that of \(s_1\), traversing each
visited layer by a path of \(A\) and splicing consecutive layers across
the matching; a layer visited twice by \(W\) is covered by two
vertex-disjoint \(A\)-paths.

\textbf{Both factors non-bipartite.} Both are then Hamilton-connected,
and no paired-cover or parity device is needed. Since \(s_0=(a_0,b_0)\)
and \(s_1=(a_1,b_1)\) are distinct, they differ in some coordinate.
Suppose \(a_0\ne a_1\); the other case is symmetric. Choose a Hamilton
path \[
a_0=x_0,x_1,\ldots,x_d=a_1
\] of \(A\). For this case it is convenient to reverse the layer
convention and view \(G\) as copies of \(B\) indexed by \(V(A)\); we
traverse these \(B\)-layers in the order \(x_0,\ldots,x_d\). Since a
non-bipartite Hamilton-connected graph has at least three vertices,
choose distinct-endpoint splice coordinates \(z_0,\ldots,z_{d-1}\) as
follows. If \(d=1\), choose \(z_0\notin\{b_0,b_1\}\). If \(d>1\), choose
\(z_0\ne b_0\), then choose \(z_i\ne z_{i-1}\) for \(1\le i<d-1\), and
finally choose \(z_{d-1}\notin\{z_{d-2},b_1\}\). Each choice excludes at
most two vertices of \(B\).

In the layer indexed by \(x_0\), take a Hamilton path of \(B\) from
\(b_0\) to \(z_0\). In each internal layer \(x_i\), \(1\le i<d\), take
one from \(z_{i-1}\) to \(z_i\), and in the last layer \(x_d\) take one
from \(z_{d-1}\) to \(b_1\). Every pair of endpoints is distinct by
construction, so Hamilton-connectedness applies in every layer. For
\(0\le i<d\), splice the two consecutive layers with the product edge
\((x_i,z_i)(x_{i+1},z_i)\). The displayed layer sequence covers each
\(B\)-layer once and its splice vertices are distinct within every
internal layer, so the result is one Hamilton \(s_0\)--\(s_1\) path.

\textbf{One factor bipartite (the substantive case).} Say \(A\) is
bipartite; the other order is symmetric. By Proposition 2.2 a nontrivial
bipartite interchange graph is balanced of even order, so \(|V(A)|=2\)
or \(|V(A)|\ge4\). Assume \(|V(A)|\ge4\) here; the case \(A=K_2\) is
treated separately after Proposition 4.1. Each color class of \(A\) then
has at least two vertices, so the Bipartite 2-DPC Lemma makes \(A\)
paired \(2\)-disjoint-path-coverable and Theorem 7.3 with \(\ell=1\)
makes it Hamilton-laceable; and \(B\) is non-bipartite and
Hamilton-connected. (Throughout this section \(n=|V(B)|\) and the
endpoints are written \(s_0=(a_0,b_0),\,s_1=(a_1,b_1)\) --- not to be
confused with the global margins.)

\begin{paperproposition}

\textbf{Proposition 4.1.} \emph{If \(A\) is balanced bipartite,
Hamilton-laceable, and paired \(2\)-disjoint-path-coverable with
\(|V(A)|\ge 4\), and \(B\) is non-bipartite and Hamilton-connected, then
\(A\,\square\,B\) is Hamilton-connected.}

\end{paperproposition}

Standard product arguments assume both factors are Hamilton-connected;
they do not apply here, because the fiber factor is bipartite, and a
Hamilton path in a balanced bipartite graph joins only vertices in
opposite color classes --- so no such factor is Hamilton-connected.
Proposition 4.1 instead manufactures Hamilton-connectedness of the
product from a factor that cannot have it, absorbing the color-parity
obstruction through the doubled layer; this is what the paired
\(2\)-disjoint-path-cover hypothesis pays for.

\emph{Proof.} The plan: choose a spanning walk of \(B\) whose parity
matches the color constraint imposed by the endpoints'
\(A\)-coordinates, then assign boundary vertices of \(A\) along the walk
so that every doubled layer receives a valid paired-cover demand. Fix
distinct \(s_0=(a_0,b_0)\) and \(s_1=(a_1,b_1)\). Because \(A\) is
bipartite, traversing a layer by a Hamilton path of \(A\) flips the
\(A\)-color of its endpoints, so threading single layers along \(W\)
forces a fixed relation between \(\chi(a_0)\) and \(\chi(a_1)\)
(throughout, \(\chi\) is applied to \(A\)-coordinates) that need not
hold. To absorb the discrepancy we \textbf{double} a layer (occasionally
two, in the same-\(B\)-layer case treated below): cover it by two
vertex-disjoint \(A\)-paths \(P,Q\) forming a paired \(2\)-disjoint path
cover of \(A\), joined by a detour into an adjacent layer. Figure 1
shows the device.

\begin{figure}
\centering
\includegraphics[width=0.93\textwidth,height=\textheight]{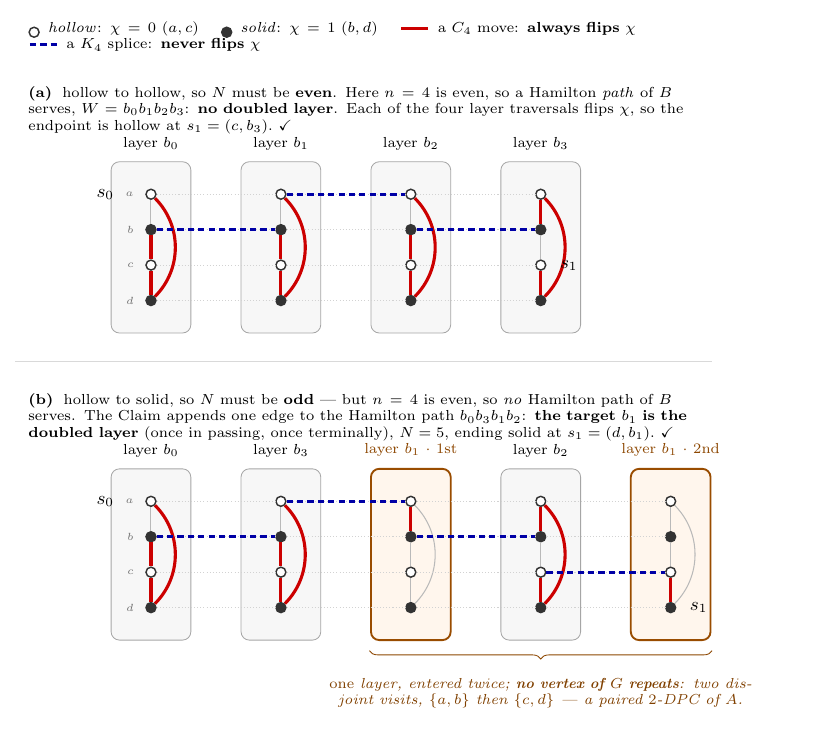}
\caption{\textbf{The doubled-layer device, as a parity fork.} The
instance is \(A=C_4\) over \(B=K_4\), with \(s_0=(a,b_0)\). Each panel
shows the layer occurrences of the selected walk, in walk order, joined
by the identity matching on the \(A\)-coordinate; edges of
\(C_4\,\square\,K_4\) not used by that presentation are suppressed, and
in panel (b) the two boxes labeled \(b_1\) are two occurrence views of
one layer, and its two visits use disjoint vertex sets. Vertex fill
gives the \(A\)-color (hollow \(\chi=0\), solid \(\chi=1\)); edge color
gives the factor supplying a move (red inside a layer, blue for a
splice). Because \(A\) is balanced bipartite, threading the layers
forces the endpoint parity \(\chi(a_1)\equiv\chi(a_0)+N\pmod 2\), where
\(N\) is the number of layer occurrences along the walk (the mechanism
is developed below). \textbf{(a)} When the demanded parity matches that
of \(n=|V(B)|=4\), a Hamilton path of \(B\) serves and
Hamilton-laceability of \(A\) suffices. \textbf{(b)} When it is
opposite, the \(N=n+1\) controlled walk constructed below re-enters the
target layer \(b_1\), doubling exactly that layer; its two visits are
the paired \(2\)-disjoint path cover the proposition hypothesizes. The
full argument is given below. When \(b_0=b_1\) the selected walk is
closed, with \(N=n+1\) or \(N=n+2\) according to the parity required; in
the \(N=n+2\) case its two doubled layers may share a second boundary,
which the proof's final case handles separately.}
\end{figure}

Write a spanning walk as \(W=(w_0,\ldots,w_{N-1})\), where \(N\) is its
number of layer occurrences, and write \(\alpha_0,\ldots,\alpha_N\) for
the \(A\)-coordinates at the boundaries of the layer pieces. Thus the
piece in occurrence \(w_i\) runs from \(\alpha_i\) to \(\alpha_{i+1}\),
with \(\alpha_0=a_0\) and \(\alpha_N=a_1\). Every single-layer Hamilton
path and both paths in every doubled-layer cover have opposite-colored
ends. Hence \[
\chi(\alpha_{i+1})=1-\chi(\alpha_i)\quad(0\le i<N),
\] and telescoping gives the single requirement \[
N\equiv\chi(a_0)+\chi(a_1)\pmod2.
\]

The doubled-layer device needs spanning walks of \(B\) of a prescribed
parity repeating no vertex more than twice. The next claim supplies
them, using only that \(B\) is Hamilton-connected (hence, on \(\ge3\)
vertices, \(2\)-connected).

\begin{paperclaim}

\textbf{Claim (controlled spanning walks).} \emph{Let \(B\) be
Hamilton-connected with \(n=|V(B)|\ge3\). For any \(s,t\in V(B)\) there
are spanning walks of \(B\) visiting every vertex at most twice and
realizing \textbf{both} parities of the occurrence count \(N\): open
walks \(s\to t\) with \(N=n\) and \(N=n+1\) when \(s\ne t\), and closed
walks based at \(s\) with \(N=n+1\) and \(N=n+2\) when \(s=t\).}

\end{paperclaim}

\emph{Proof of Claim.} The needed neighbors can be read off Hamilton
paths directly. If \(s\ne t\): a Hamilton path \(s\to t\) gives \(N=n\);
its \emph{penultimate} vertex \(w\) is a neighbor of \(t\) with
\(w\ne s\) (a Hamilton path on \(n\ge3\) vertices visits each vertex
once, and \(w\) is not the first), so a Hamilton path \(s\to w\)
followed by the edge \(wt\) gives \(N=n+1\) with \(t\) the only repeated
vertex. If \(s=t\): the second vertex \(u\) of any Hamilton path out of
\(s\) is a neighbor of \(s\), and a Hamilton path \(s\to u\) closed by
the edge \(us\) gives \(N=n+1\); taking a neighbor \(v\) of \(s\) and
the penultimate vertex \(u'\) of a Hamilton path \(s\to v\) (a neighbor
of \(v\) with \(u'\ne s\)), a Hamilton path \(s\to u'\) followed by
\(u'\,v\,s\) gives \(N=n+2\), with \(s\) (the required return) and \(v\)
the only repeated vertices. \(\square\)

The four Claim constructions have, respectively, zero, one, one, and two
doubled layers. In the open \(N=n+1\) walk, the earlier occurrence of
\(b_1\) is internal in the fresh Hamilton path, at some position
\(1\le j\le n-2\), so its boundary pair is disjoint from the final pair.
In the closed \(N=n+1\) walk, the two occurrences of \(b_0\) are first
and last and their boundary pairs are disjoint because \(n\ge3\). If
\(b_0=b_1\), then \(a_0\ne a_1\), since the two product endpoints are
distinct. The remaining two-doubled-layer incidence is analyzed below.

Apply the claim with \((s,t)=(b_0,b_1)\). The two walks have opposite
\(N\)-parities, so one satisfies the displayed requirement; take it as
\(W\). Its repeated vertices are the doubled layers: at most \(b_1\)
when \(b_0\ne b_1\), and at most \(b_0\) together with one other layer
when \(b_0=b_1\). Every singly visited layer will receive a Hamilton
path of \(A\). We now make the global boundary choice that supplies
valid demands in all doubled layers.

\begin{paperlemma}

\textbf{Lemma 4.2 (boundary assignment).} \emph{For each controlled walk
in the Claim, there is an assignment \(\alpha_0,\ldots,\alpha_N\) with
\(\alpha_0=a_0\), \(\alpha_N=a_1\), and
\(\chi(\alpha_{i+1})\ne\chi(\alpha_i)\) for every \(i\), such that the
four terminals belonging to each doubled layer are pairwise distinct.
Thus its two occurrence pairs form a pairwise-opposite demand for a
paired \(2\)-cover of \(A\).}

\end{paperlemma}

\emph{Proof.} In the one-doubled-layer constructions the doubled layer
contains at most one of the endpoint pins unless \(b_0=b_1\); in that
closed-walk case it contains both, and then \(a_0\ne a_1\) because
\(s_0\ne s_1\).

The colors of all boundaries are forced by \(\chi(a_0)\) and their
indices. With no doubled layer, choose any vertex of the required color
at each free boundary. With one doubled layer, its four terminals
contain at most the two distinct pins \(a_0,a_1\). Their forced colors
occur twice each. If the two pins have the same color, they already
occupy the two terminals of that color and we choose two distinct
vertices of the other color. If the pins have opposite colors, choose in
each color a second vertex distinct from its pin. With one pin or none,
fewer terminals are forced, so the same choice succeeds with more room
to spare. Each color class has at least two vertices, so all choices
exist.

It remains to treat the closed walk with \(N=n+2\), where the two
doubled layers may have two common boundaries. Use the notation of the
Claim. Let \[
x_0=b_0,x_1,\ldots,x_{n-1}=u'
\] be the fresh Hamilton path from \(b_0\) to \(u'\), and let \(v=x_k\).
Since \(v\ne b_0,u'\), we have \(1\le k\le n-2\), and \[
W=(x_0,x_1,\ldots,x_{n-1},v,b_0).
\] The doubled \(b_0\)-layer has terminal pairs \[
(\alpha_0,\alpha_1),\qquad(\alpha_{n+1},\alpha_{n+2}),
\] and the doubled \(v\)-layer has \[
(\alpha_k,\alpha_{k+1}),\qquad(\alpha_n,\alpha_{n+1}).
\] The layers always share \(\alpha_{n+1}\). If \(k=1\), equivalently if
the fresh Hamilton path begins \(b_0,v,\ldots\), they also share
\(\alpha_1\). This is the two-shared-boundary configuration; for
\(B=K_3\) it is forced. When \(k=n-2\), one has \(k+1=n-1<n\), so
\(\alpha_{k+1}\) is not the boundary occurrence \(\alpha_n\) and no
further shared boundary is created. For \(n=3\), one has \(k=1=n-2\),
and the two-shared-boundary configuration just identified is the unique
detour shape.

Thus the fixed terminals are \(\alpha_0=a_0\) and \(\alpha_{n+2}=a_1\);
the two demands always share \(\alpha_{n+1}\), and share \(\alpha_1\)
only when \(k=1\), every other displayed terminal belonging to just one
of them. We therefore choose the shared boundaries first.

The choices are made in the following order, each boundary taken in the
color its index forces and avoiding at most one already-chosen vertex of
that color: \[\begin{array}{l|l|l}
\text{boundary} & \text{forced color} & \text{must avoid}\\\hline
\alpha_{n+1} & \chi(a_0)+n+1 & \alpha_0=a_0\\
\alpha_1 & \chi(a_0)+1 & \text{whichever of }\alpha_{n+1},\alpha_{n+2}\text{ shares its color}\\
\alpha_n & \chi(a_0)+n & \alpha_1,\ \text{when their colors agree}\\
\alpha_k,\ \alpha_{k+1} & \chi(a_0)+k,\ \chi(a_0)+k+1 & \text{the same-colored member of }\{\alpha_n,\alpha_{n+1}\}
\end{array}\] Colors are read modulo \(2\); \(\alpha_0=a_0\) and
\(\alpha_{n+2}=a_1\) are pinned --- fixed in advance, rather than chosen
by the construction --- and when \(k=1\) the entry \(\alpha_k=\alpha_1\)
has already been chosen. In words:

Choose the relevant boundaries in the order used by the formal
construction. First choose \(\alpha_{n+1}\) in its forced color,
distinct from \(a_0=\alpha_0\). Next choose \(\alpha_1\) in its forced
color. The vertices \(a_1=\alpha_{n+2}\) and \(\alpha_{n+1}\) have
opposite colors, so exactly one of them has the color of \(\alpha_1\);
avoid that one. The four terminals of the \(b_0\)-layer are now
distinct. Choose \(\alpha_n\) in its forced color, avoiding \(\alpha_1\)
when their colors agree.

If \(k=1\), set \(\alpha_k=\alpha_1\) and choose \(\alpha_{k+1}\) in its
forced color, avoiding the one member of \(\{\alpha_n,\alpha_{n+1}\}\)
having that color. The four terminals of the \(v\)-layer are then
distinct. If \(k>1\), choose \(\alpha_k\) and \(\alpha_{k+1}\) in their
forced, opposite colors; each avoids the unique same-colored member of
\(\{\alpha_n,\alpha_{n+1}\}\). Choose all remaining boundaries
arbitrarily in their forced colors.

This order also proves the extremal case in which each color class of
\(A\) has exactly two vertices. At every choice, at most one previously
chosen vertex of the required color is excluded; vertices of the
opposite color are automatically distinct. Hence a choice remains in a
two-element color class. Both the one- and two-shared-boundary
configurations are covered. \(\square\)

With the boundaries fixed, each singly-visited layer is traversed by a
Hamilton path of \(A\) between its opposite-colored boundaries
(Hamilton-laceability), and each doubled layer --- including the layer
carrying both pins \(a_0\ne a_1\) in the same-\(B\)-layer case --- by a
paired \(2\)-cover, whose four terminals form a pairwise-opposite demand
by Lemma 4.2, so the paired-cover hypothesis on \(A\) supplies it.
Splicing all the layer pieces across the matching along \(W\) yields the
Hamilton \(s_0\)--\(s_1\) path. \(\square\)

\textbf{The degenerate factor \(A=K_2\).} When \(|V(A)|=2\) the
doubled-layer device is unavailable (a paired \(2\)-cover needs four
distinct terminals), and \(G=K_2\,\square\,B\) is the prism over \(B\):
two copies \(B^0,B^1\) joined by the rungs \(v^0v^1\). We prove it
Hamilton-connected directly, using only that \(B\) (Hamilton-connected
by the induction hypothesis, with \(n=|V(B)|\ge3\)) has a Hamilton path
between any two distinct vertices. Let the endpoints be \((i,u)\) and
\((j,w)\).

\begin{itemize}
\tightlist
\item
  \emph{Different copies} (\(i\ne j\), say \((0,u)\to(1,w)\)): pick
  \(x\notin\{u,w\}\) (possible as \(n\ge3\); if \(u=w\) any \(x\ne u\)).
  A Hamilton path \(u\to x\) of \(B^0\), the rung \(x^0x^1\), and a
  Hamilton path \(x\to w\) of \(B^1\) concatenate into a Hamilton path
  of the prism.
\item
  \emph{Same copy} (\(i=j\), say \((0,u)\to(0,w)\) with \(u\ne w\)):
  take a Hamilton path \(u=h_0,h_1,\dots,h_{n-1}=w\) of \(B\). Then
  \[\begin{aligned}&(0,h_0)\to(1,h_0)\to[\,\text{a fresh Hamilton path }h_0\to h_1\text{ of }B^1\,]\\&\qquad\to(1,h_1)\to(0,h_1)\to(0,h_2)\to\cdots\to(0,h_{n-1})\end{aligned}\]
  is a Hamilton path of the prism: \(B^1\) is covered by the interior
  Hamilton path \(h_0\to h_1\), and \(B^0\) by \((0,h_0)\) together with
  the tail \((0,h_1),\dots,(0,h_{n-1})\).
\end{itemize}

Hence every pair of prism vertices is joined by a Hamilton path, so
\(K_2\,\square\,B\) is Hamilton-connected. \(\square\)

\begin{center}\rule{0.5\linewidth}{0.5pt}\end{center}

\hypertarget{large-non-bipartite-prime-classes}{%
\subsection{5. Large non-bipartite prime
classes}\label{large-non-bipartite-prime-classes}}

This section handles case (iii), the largest branch and the only
genuinely new construction: the indecomposable non-bipartite prime
classes. It applies the induction hypothesis to the fibers of a
separating line and uses the Naddef--Pulleyblank dichotomy to obtain
Hamilton-connectedness of the non-bipartite matroid base-exchange
quotient.

Let \(G=G(R,S)\) be a nonempty class normalized by deleting inactive
lines, and assume that it is invariant-free (hence prime/indecomposable
by §2.2), non-bipartite, has \(m,n\ge3\) active rows and columns, and
has \(|V(G)|>6\). Fix distinct vertices \(a,b\in V(G)\); these are the
\textbf{pinned} endpoints, the two vertices the Hamilton path under
construction must join, and a vertex or a color is called pinned when it
is fixed in advance in this way rather than chosen by the construction.
For a line \(L\), write \(F_p\) for the fiber obtained by fixing \(L\)
to the realizable pattern \(p\) (§2.1), and write \(Q_L\) for its
quotient. Four notions recur:

\begin{itemize}
\tightlist
\item
  a line \(L\) is \textbf{separating} (for \(a,b\)) if \(a\) and \(b\)
  have different patterns on \(L\);
\item
  the \textbf{buffer line} is the separating line produced by Lemma 5.9,
  one of whose fibers --- the \textbf{buffer} --- is non-bipartite; it
  is the line we build the final path around. After transposing if
  necessary we take the buffer line to be a row, the \textbf{pivot row},
  and every statement below has the transposed column form; a row other
  than the pivot row is a \textbf{non-pivot row};
\item
  for adjacent patterns \(p,q\) of \(L\), a \textbf{crossing edge}
  (equivalently, an \textbf{interface edge}) is an edge of \(G\) with
  one endpoint in \(F_p\) and the other in \(F_q\); the endpoints of
  such edges form the \textbf{interface} between the two fibers;
\item
  a quotient edge \(pq\) is \textbf{used} if it is one of the edges of
  the Hamilton path of \(Q_L\) chosen in step (3) below, so that the
  construction must cross between \(F_p\) and \(F_q\) along it.
\end{itemize}

Section 6 is logically prior to this one: its base results are among the
strong-induction inputs used below, and are presented afterward only to
keep the pivot construction together. Every lemma in this section is
stated within this regime (active, invariant-free, indecomposable,
non-bipartite, \(m,n\ge3\)). The order bound \(|V(G)|>6\) selects the §5
branch rather than the §6 branch; it is not used in the local proofs of
Lemmas 5.8--5.10.

The whole construction is four steps, each licensed by the lemma named
beside it; the reader may hold this skeleton while the lemmas supply the
guarantees; it is restated as an explicit algorithm and run on a worked
example at the end of the section. Given distinct endpoints \(a,b\): (1)
choose a separating \textbf{buffer line} \(L\) with a non-bipartite
fiber (Lemma 5.9); (2) its quotient \(Q_L\) is a non-bipartite matroid
base-exchange graph, hence Hamilton-connected (Lemmas 5.3, 5.4, 5.10,
Corollary 5.5); (3) take a Hamilton path of \(Q_L\) from the fiber of
\(a\) to the fiber of \(b\); (4) thread the fibers along that path,
traversing each by a recursively obtained Hamilton path and splicing
consecutive fibers on interface edges (Lemmas 5.6, 5.7, Proposition
5.11). The buffer fiber absorbs the parity that a chain of forced fibers
could not.

Figure 2 shows the whole construction on a small class, and the reader
may find it useful to hold while the lemmas arrive; it is the class
worked in full at the end of the section.

\begin{figure}
\centering
\includegraphics[width=1\textwidth,height=\textheight]{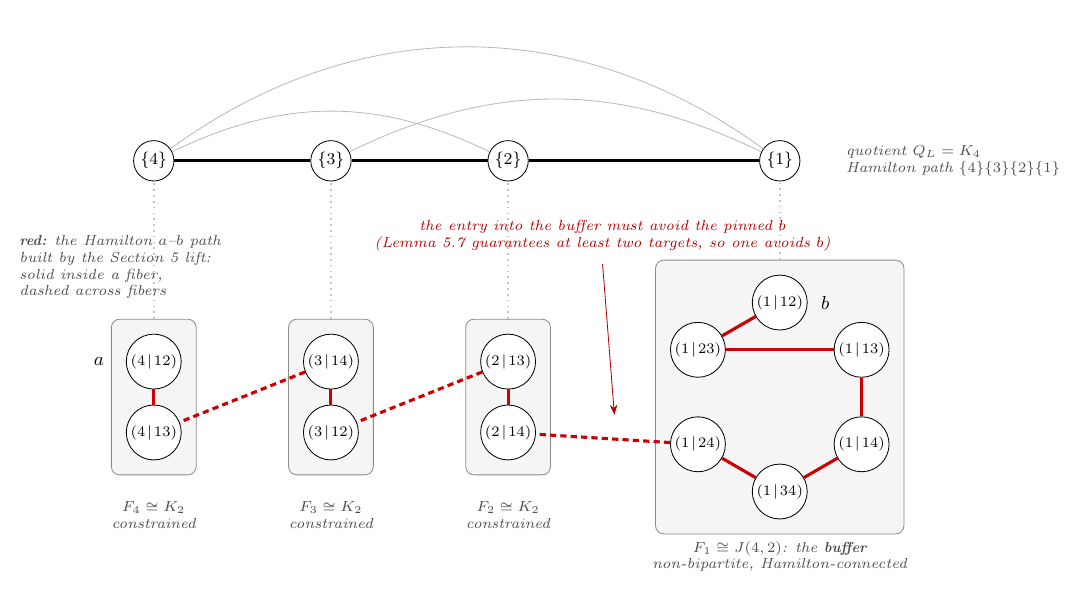}
\caption{\textbf{The construction of this section, on the class
\(R=(2,2,1)\), \(S=(2,1,1,1)\).} Here \((c\mid xy)\) denotes the matrix
whose row-\(3\) entry lies in column \(c\) and whose row-\(1\) support
is \(\{x,y\}\); the margins then determine row \(2\). All twelve
matrices of \(\mathcal A(R,S)\) are drawn; of the thirty-three edges of
\(G(R,S)\) only the eleven the lift uses are shown, the other twenty-two
being omitted for legibility. Fixing the buffer line \(L=\) row \(3\)
partitions the class into four fibers, which sit below their vertices in
the quotient \(Q_L=K_4\) (dotted). Three fibers are \(K_2\)'s, hence
bipartite, so a traversal's exit is forced (such fibers are called
\emph{constrained} below); the fourth, \(F_1\cong J(4,2)\), is the
non-bipartite buffer. A Hamilton path of \(Q_L\) (heavy, above) is
lifted to a Hamilton \(a\)--\(b\) path of \(G\) (red): each fiber is
traversed by a Hamilton path obtained from the induction hypothesis, and
consecutive fibers are spliced on crossing edges (red dashed). Of the
eighteen crossing edges of \(G(R,S)\) the lift uses three. For these
endpoints and this quotient path the first two splices are forced once
the constrained-fiber exits are fixed, and what Lemmas 5.6 and 5.7
guarantee there is that a usable crossing exists at all. The last
crossing is the delicate one: it must enter the buffer without landing
on the pinned endpoint \(b\). Lemma 5.7 supplies two targets from the
forced exit \((2\mid14)\), namely \((1\mid14)\) and \((1\mid24)\), and
the lift takes \((1\mid24)\).}
\end{figure}

Throughout this section every fiber is immediately put into
\textbf{canonical normal form}: delete inactive lines, then split the
invariant positions into their Cartesian factors (§2.2). We call the
result the fiber's \textbf{canonical core}. The induction hypothesis
applies to any strictly smaller interchange graph --- in particular to
the resulting invariant-free factors, and to whole fiber cores where
convenient. In particular, the induction hypothesis includes every fiber
of order at most six and every Johnson fiber. The smaller instance
follows the ordinary §3 dispatch: a bipartite fiber is handled outright
by Lemma 1.2, a decomposable non-bipartite fiber by §4, a Johnson fiber
by the §6 Johnson branch, and a remaining prime many-line fiber of order
at most six by the §6 finite branch. Thus an order-six fiber inside a
larger §5 class is already covered by the strong induction, not a new §5
case.

\hypertarget{normalizing-fibers-and-quotients-lemmas-5.15.4-and-corollary-5.5}{%
\subsubsection{Normalizing Fibers and Quotients (Lemmas 5.1--5.4 and
Corollary
5.5)}\label{normalizing-fibers-and-quotients-lemmas-5.15.4-and-corollary-5.5}}

\hypertarget{canonical-fibers}{%
\paragraph{Canonical fibers}\label{canonical-fibers}}

\begin{paperlemma}

\textbf{Lemma 5.1 (canonical fiber-core renormalization).}\\
Let \(L\) be a separating line and let \(F_p\) be a fiber of \(L\).
After deleting \(L\) and canonicalizing \(F_p\) by deleting inactive
lines and splitting invariant positions, every nontrivial factor is a
strictly smaller interchange graph.

\end{paperlemma}

\emph{Proof.} Fixing \(L\) to \(p\) gives a proper subclass, since \(L\)
is separating and hence at least two line patterns occur in \(G\).
Deleting \(L\) leaves a complete-margin \(0\)-\(1\) matrix class on
fewer lines. Any inactive line in this class is forced in every matrix
of the fiber, so deleting it is graph-isomorphic. Any remaining
invariant position gives the Cartesian factorization of §2.2: the
non-invariant cells split into diagonal blocks and the off-block
rectangles are forced. Thus the fiber graph is a Cartesian product of
the interchange graphs of those variable blocks.

Each block is strictly smaller than \(G\): the fiber has fewer vertices
than \(G\), and every nontrivial factor has at most the fiber order.
\(\square\)

\hypertarget{pivot-freedom}{%
\paragraph{Pivot freedom}\label{pivot-freedom}}

Pivot-freedom is used twice below: in Lemma 5.9, to vary a chosen cell
outside the \emph{wide pair} (a pair of rows whose supports are
incomparable and differ in at least three columns; defined with Lemma
5.8), and in Lemma 5.10, to make every line matroid \textbf{loopless}
and \textbf{coloopless} --- that is, to rule out an element lying in no
basis and an element lying in every basis.

\begin{paperlemma}

\textbf{Lemma 5.2 (pivot-freedom).}\\
In a nonempty active prime block, every cell --- every row--column
position \((i,j)\) --- is \textbf{variable}, that is, takes the value
\(0\) in some matrix of the class and the value \(1\) in another.
Equivalently, the bipartite graph on the rows and columns whose edges
are the variable cells is the complete bipartite graph \(K_{m,n}\).

\end{paperlemma}

\emph{Proof.} This is the theorem of Brualdi and Manber \citep[Thm.
9]{brualdiPrimeInterchangeGraphs1983}: an active prime class has no
invariant cell --- every row-column position varies in some realization.
(By the prime criterion of §2.2 the block's primeness is equivalent to
its being invariant-free; the content here is that primeness forces
\emph{every} cell to vary.) \(\square\)

\emph{Remark.} One can also derive this from the theory of forced
entries (Brualdi \citep[Ch. 3]{brualdiCombinatorialMatrixClasses2006}):
a forced cell is exactly a cell separated from the variable part by a
tight row-column cut; in an active class all inactive-line forced cells
have been deleted, and in an invariant-free block no nontrivial forced
off-diagonal rectangle remains --- otherwise the block would split as a
nontrivial Cartesian product.

\hypertarget{line-quotients-as-matroid-base-graphs}{%
\paragraph{Line quotients as matroid base
graphs}\label{line-quotients-as-matroid-base-graphs}}

Let \(L\) be a row; the column case is obtained by transposition. Let
\(k=r_L\), the prescribed sum of the row \(L\), let \(R'\) be the
row-sum vector after deleting \(L\), and let \(S=(s_1,\ldots,s_n)\).
Write \([n]=\{1,\ldots,n\}\) for the set of column indices, so that the
pattern of \(L\) (§2.1) --- the set of columns in which \(L\) carries a
\(1\) --- is a \(k\)-subset of \([n]\). A pattern \(p\subseteq[n]\),
\(|p|=k\), is realizable exactly when the residual column vector
\(S-\mathbf 1_p\) is realizable with row sums \(R'\). Here
\(\mathbf 1_p\) is the indicator vector of \(p\),
\(e_i=\mathbf 1_{\{i\}}\) is the \(i\)th unit vector, and \(p_i\)
denotes the \(i\)th coordinate of \(\mathbf 1_p\), so that \(p_i=1\)
means \(i\in p\). For a pattern \(p\) we write \(p-i+j\) for
\((p\setminus\{i\})\cup\{j\}\), and reserve \(e_i\) for identities
between vectors.

The criterion is easiest to read through two set functions. For a set
\(X\) of columns, \(h(X)\) is the largest number of \(1\)s the deleted
rows can place in \(X\), and \(\ell(X)\) is therefore the shortfall that
the row \(L\) must itself cover in \(X\). Formally, for
\(X\subseteq[n]\) put \[
h(X)=\sum_i \min(r'_i,|X|),\qquad \ell(X)=\sum_{j\in X}s_j-h(X).
\] By the Gale--Ryser theorem
\citep{galeTheoremFlowsNetworks1957, ryserCombinatorialPropertiesMatrices1957}
(see also Brualdi \citep[Ch. 3]{brualdiCombinatorialMatrixClasses2006}),
\(p\) is realizable if and only if \[
|p|=k,\qquad |p\cap X|\ge \ell(X)\quad\text{for every }X\subseteq[n].
\] That is: \(p\) has the right size, and on every set of columns it
meets the shortfall.

Two properties of these functions carry the whole of Lemma 5.3. First,
\(h(X)\) depends only on \(|X|\) and is submodular, since
\(q\mapsto\sum_i\min(r'_i,q)\) is discrete concave; hence \(\ell\) is
supermodular. Second, two ambient facts let the display be stated as it
is: \(\sum_j s_j=k+\sum_i r'_i\), because the class is nonempty, and
every column is active in the §5 regime, so \(s_j\ge1\) for every \(j\)
and \(S-\mathbf 1_p\) is nonnegative for every pattern \(p\).

\begin{paperlemma}

\textbf{Lemma 5.3 (single-line patterns are matroid bases).}\\
The realizable patterns of \(L\) are the bases of a matroid \(M_L\).

\end{paperlemma}

\emph{Proof.} Let \(\mathcal B\) be the family of \(k\)-subsets
satisfying the Gale--Ryser inequalities above. We verify the
basis-exchange axiom. Let \(B,B'\in\mathcal B\) and
\(e\in B\setminus B'\). Put \(B_0=B\setminus\{e\}\). A set \(X\) is
deficient for \(B_0\) if \(|B_0\cap X|<\ell(X)\). Such an \(X\) must
contain \(e\), and since \(B\) was feasible, necessarily
\(|B_0\cap X|=\ell(X)-1\).

Deficient sets are closed under intersection. If \(X,Y\) are deficient,
then using modularity of cardinality and supermodularity of \(\ell\),
\[\begin{aligned}
|B_0\cap(X\cap Y)|
&=|B_0\cap X|+|B_0\cap Y|-|B_0\cap(X\cup Y)|\\
&\le \ell(X)+\ell(Y)-2-(\ell(X\cup Y)-1)
\le \ell(X\cap Y)-1
\end{aligned}\] (the middle step also uses
\(|B_0\cap(X\cup Y)|\ge\ell(X\cup Y)-1\), which holds for \emph{every}
set: \(B\) satisfies every inequality, and deleting the single element
\(e\) lowers each left-hand side by at most one). The family of
deficient sets is nonempty: the ground set \([n]\) is itself deficient,
since \(h([n])=\sum_i r'_i\) gives \(\ell([n])=k\) while
\(|B_0\cap[n]|=|B_0|=k-1\). Thus the intersection \(Z\) of all deficient
sets (a finite, intersection-closed, nonempty family) is deficient.
Since \(B'\) is feasible, \[
|B'\cap Z|\ge \ell(Z)=|B_0\cap Z|+1.
\] As \(e\notin B'\), this gives an element \(f\in B'\setminus B\) with
\(f\in Z\). Because \(Z\) is contained in every deficient set, \(f\)
lies in every deficient set, so adjoining it raises the left side of
each deficient inequality from \(\ell(X)-1\) to \(\ell(X)\); the
non-deficient inequalities were already satisfied by \(B_0\). Hence
\(B-e+f\in\mathcal B\). \(\square\)

\emph{Remark.} The matroid property of single-line patterns is classical
in spirit --- by Gale--Ryser the realizable patterns are the integer
points of a polymatroid slice, and such families form matroid bases
(Edmonds \citep{edmondsSubmodularFunctionsMatroids1970}) --- and we keep
the short self-contained proof above. The lemma also admits an
independent proof by a minimal-pair argument: among realizations
\(M\in F_B\), \(N\in F_{B'}\) choose a pair minimizing the number of
differing cells; two local interchange arguments then force a closure
property on the difference pattern, and a counting identity over the
difference cells produces the exchange as a single explicit interchange.
(The same minimal-pair technique proves the unit-transfer Lemma 5.4
below, whose proof we give in that form.)

Given two columns \(i,j\) and a matrix \(X\), we say that a row of \(X\)
has \textbf{type} \((\varepsilon,\varepsilon')\) \textbf{in columns}
\((i,j)\) when its entries in those two columns are \(\varepsilon\) and
\(\varepsilon'\) respectively.

\begin{paperlemma}

\textbf{Lemma 5.4 (unit-transfer).} \emph{If \(p\) and \(q=p-i+j\) are
both realizable patterns of \(L\) (so \(i\in p\), \(j\notin p\)), then
some matrix in \(F_p\) has a non-pivot row of type \((0,1)\) in columns
\((i,j)\).}

\end{paperlemma}

\emph{Proof.} Write \(c=S-\mathbf 1_p\) for the residual column-sum
vector of the non-pivot rows when \(L=p\), and
\(c'=S-\mathbf 1_q=c+e_i-e_j\); the residual classes
\(\mathcal A'=\mathcal A(R',c)\) and \(\mathcal B'=\mathcal A(R',c')\)
are both nonempty, by the realizability of \(p\) and of \(q\). Suppose,
for contradiction, that no matrix of \(\mathcal A'\) has a row of type
\((0,1)\) in columns \((i,j)\): then in every member of \(\mathcal A'\),
every row with a \(1\) in column \(j\) also has a \(1\) in column \(i\).

Choose a pair \((M,N)\in\mathcal A'\times\mathcal B'\) minimizing the
number of cells at which \(M\) and \(N\) differ. Call a cell a
\emph{P-cell} if \(M\) has a \(1\) and \(N\) a \(0\) there, and a
\emph{Q-cell} in the opposite case. In each column \(k\), the count of
P-cells minus Q-cells equals \(c_k-c'_k\): zero for \(k\ne i,j\), \(+1\)
for \(k=j\) (and \(-1\) for \(k=i\), which we shall not need). Hence the
row set \[
\mathcal S=\{a:\ (a,j)\text{ is a P-cell}\}
\] is nonempty, and every \(a\in\mathcal S\) has \(M_{aj}=1\), so also
\(M_{ai}=1\) by our assumption. Let \(\mathcal T\) be the set of columns
containing a Q-cell of some row of \(\mathcal S\). Then
\(j\notin\mathcal T\) and \(i\notin\mathcal T\) (rows of \(\mathcal S\)
carry \(1\)s of \(M\) in both columns).

\emph{Closure: for \(k\in\mathcal T\), every P-cell of column \(k\) lies
in a row of \(\mathcal S\).} Let \((b,k)\) be a P-cell and pick
\(a\in\mathcal S\) with a Q-cell at \((a,k)\); then \(a\ne b\) and
\(k\ne j\) (row \(a\) has \(M_{aj}=1\) but \(M_{ak}=0\)). The known
cells so far: \[\begin{array}{c|ccc}
 & (a,j) & (a,k) & (b,k)\\\hline
M & 1 & 0 & 1\\
N & 0 & 1 & 0
\end{array}\] If \(M_{bj}=0\), the block on rows \(\{a,b\}\) and columns
\(\{j,k\}\) of \(M\) is switchable (\(M_{aj}=M_{bk}=1\),
\(M_{ak}=M_{bj}=0\)); the switch stays in \(\mathcal A'\) and flips four
cells of which at least three, namely \((a,j)\), \((a,k)\), \((b,k)\),
were differences from \(N\), so it strictly decreases the difference
count, contradicting minimality. Hence \(M_{bj}=1\).

If \(N_{bj}=1\), the block on rows \(\{b,a\}\) and columns \(\{j,k\}\)
of \(N\) is switchable (\(N_{bj}=N_{ak}=1\), \(N_{bk}=N_{aj}=0\)); the
switch stays in \(\mathcal B'\) and again flips at least three current
differences, namely \((b,k)\), \((a,j)\), \((a,k)\), contradicting
minimality. Hence \(M_{bj}=1\) and \(N_{bj}=0\),
i.e.~\(b\in\mathcal S\).

\emph{Counting.} Fix \(a\in\mathcal S\). Rows of \(M\) and \(N\) have
equal sums, so row \(a\) has exactly as many P-cells as Q-cells. All its
Q-cells lie in \(\mathcal T\)-columns (by the definition of
\(\mathcal T\)), while among its P-cells the cell \((a,j)\) lies outside
\(\mathcal T\); hence its P-cells in \(\mathcal T\)-columns number at
most its Q-cells in \(\mathcal T\)-columns minus one. Summing over
\(\mathcal S\) and recounting by columns, \[
\sum_{k\in\mathcal T}P^{\mathcal S}_k+|\mathcal S|\ \le\ \sum_{k\in\mathcal T}Q^{\mathcal S}_k,
\] where \(P^{\mathcal S}_k,Q^{\mathcal S}_k\) count the P- and Q-cells
of column \(k\) in rows of \(\mathcal S\), and \(P_k,Q_k\) denote the
corresponding counts over all rows. But each \(k\in\mathcal T\) has
\(k\ne i,j\), so column \(k\) has equally many P- and Q-cells in total;
the closure property gives \(P^{\mathcal S}_k=P_k\), and trivially
\(Q^{\mathcal S}_k\le Q_k=P_k\). Hence the right side is at most the
left sum without the \(|\mathcal S|\) term, forcing
\(|\mathcal S|\le 0\) --- contradicting \(\mathcal S\ne\emptyset\).

Therefore some matrix of \(F_p\) has a non-pivot row of type \((0,1)\)
in columns \((i,j)\); switching that row against \(L\) on columns
\(\{i,j\}\) realizes the exchange \(p\to q\) as an actual crossing edge
between \(F_p\) and \(F_q\). \(\square\)

\emph{Remark.} Alternatively, the lemma follows from Gale--Ryser
tightness (Brualdi \citep[Ch.
3]{brualdiCombinatorialMatrixClasses2006}): column \(j\) is contained in
column \(i\) in every member of a nonempty class exactly when some
column set \(X\) with \(i\in X\), \(j\notin X\) is Gale--Ryser-tight
(its inequality above holds with equality), equivalently when the
transferred margin vector \(c+e_i-e_j\) is infeasible --- a criterion
derived there from the Gale--Ryser theorem; applying it with the
residual margins of \(p\) gives the row sought.

By Lemma 5.4 every basis exchange of \(M_L\) is realized by a crossing
edge, and conversely every crossing interchange changes the line pattern
by exactly one basis exchange. Therefore \(Q_L\) \textbf{equals} the
basis-exchange graph of \(M_L\), i.e.~the \(1\)-skeleton of the matroid
base polytope \(\operatorname{conv}\{\mathbf 1_B:B\in\mathcal B(M_L)\}\)
(two bases are polytope-adjacent exactly when they differ by a single
exchange --- a criterion of Edmonds, first in print in Hausmann--Korte
\citep[Theorem 4.3]{hausmannColouringCriteriaAdjacency1978}; see also
\citep[p.~309]{NP1981}).

\begin{papercorollary}

\textbf{Corollary 5.5.}\\
If \(Q_L\) is non-bipartite, then \(Q_L\) is Hamilton-connected.

\end{papercorollary}

\emph{Proof.} By Lemmas 5.3--5.4, \(Q_L\) is the \(1\)-skeleton of a
matroid base polytope. By Naddef--Pulleyblank \citep[Theorem 3.3.1 and
Corollary 3.3.2]{NP1981}, the graph of a matroid base polytope is either
a hypercube or Hamilton-connected. Since \(Q_L\) is non-bipartite, the
hypercube alternative is impossible. \(\square\)

\hypertarget{interfaces-between-fibers-lemmas-5.65.7}{%
\subsubsection{Interfaces Between Fibers (Lemmas
5.6--5.7)}\label{interfaces-between-fibers-lemmas-5.65.7}}

Call a fiber \textbf{constrained} if its canonical core is bipartite of
order at least \(2\). By the bipartite classification it is balanced,
and by the induction hypothesis it is Hamilton-laceable. A non-bipartite
fiber is Hamilton-connected by induction; a singleton imposes no
internal parity condition. In summary:

\begin{longtable}[]{@{}
  >{\raggedright\arraybackslash}p{(\columnwidth - 4\tabcolsep) * \real{0.1531}}
  >{\raggedright\arraybackslash}p{(\columnwidth - 4\tabcolsep) * \real{0.4184}}
  >{\raggedright\arraybackslash}p{(\columnwidth - 4\tabcolsep) * \real{0.4286}}@{}}
\toprule\noalign{}
\begin{minipage}[b]{\linewidth}\raggedright
fiber type
\end{minipage} & \begin{minipage}[b]{\linewidth}\raggedright
internal traversal
\end{minipage} & \begin{minipage}[b]{\linewidth}\raggedright
interface freedom
\end{minipage} \\
\midrule\noalign{}
\endhead
\bottomrule\noalign{}
\endlastfoot
constrained & Hamilton-laceable; exit color forced & whole fiber (5.6);
two reachable targets (5.7) \\
non-bipartite & Hamilton-connected (induction) & \(\ge2\) interface
vertices (5.6); two-sided avoid (5.7) \\
singleton & trivial & no parity constraint; freedom supplied on the far
side (5.6) \\
\end{longtable}

\addtocounter{table}{-1}

(In the worked example closing the section, the buffer fiber is the
octahedron \(J(4,2)\) and the three other fibers are single edges
\(K_2\) --- the reader may find it useful to glance at that instance
while reading Lemmas 5.6--5.7.) These two lemmas are consumed only in
the gluing of Proposition 5.11; a reader may take their statements on a
first pass and return to the proofs alongside the gluing.

\begin{paperlemma}

\textbf{Lemma 5.6 (interface on a used edge).} \emph{Let \(p\) and
\(q=p-i+j\) be adjacent patterns of a pivot row \(L\) (so \(p_i=1\),
\(p_j=0\)), and let \(I_{p\to q}\) be the set of vertices of \(F_p\)
incident to a crossing edge into \(F_q\). If \(F_p\) is bipartite, then
\(I_{p\to q}=V(F_p)\). If \(F_p\) is non-bipartite and the edge is used,
then \(|I_{p\to q}|\ge 2\).}

\end{paperlemma}

\emph{Proof.} Write \(\operatorname{supp}_j(M)=\{r:M_{rj}=1\}\) for the
set of rows in which column \(j\) of \(M\) carries a \(1\). A crossing
switch \(p\to q\) at \(M\in F_p\) exists exactly when some non-pivot row
of \(M\) has type \((0,1)\) in columns \((i,j)\); so the non-interface
is the nesting locus \[
N=\{M\in F_p:\operatorname{supp}_j(M)\subseteq \operatorname{supp}_i(M)\},\qquad I_{p\to q}=V(F_p)\setminus N.
\] Suppose \(N\) is neither empty nor all of \(F_p\). Since fibers are
connected, there is an internal edge \(M M'\) with \(M\in N\) and
\(M'\notin N\). The switch \(M\to M'\) uses exactly one of the columns
\(i,j\): using neither would leave the nesting unchanged, and using both
would exhibit a \((0,1)\)-row already in \(M\). Say it uses columns
\(i,c\); then \(M\) has two rows \[
u:(i,j,c)=(1,1,0),\qquad v:(i,j,c)=(0,0,1),
\] and the alternate switch on columns \(j,c\) produces a third vertex
\(M''\notin N\) with \(M''\ne M'\). If the boundary switch instead uses
columns \(j,c\), the alternate switch on columns \(i,c\) produces the
same three-vertex configuration.

If \(F_p\) is bipartite, then \(M,M',M''\) differ pairwise by single
switches --- \(M'\) from \(M\) on \(\{i,c\}\), \(M''\) from \(M\) on
\(\{j,c\}\), and \(M'\) from \(M''\) on \(\{i,j\}\) --- a triangle,
contradicting bipartiteness. Hence \(N\) is empty or all of \(F_p\).
Since \(p,q\) are adjacent, Lemma 5.4 realizes the quotient edge, so
\(I_{p\to q}\ne\emptyset\) and \(N\ne F_p\); thus \(N=\emptyset\) and
\(I_{p\to q}=V(F_p)\).

If \(F_p\) is non-bipartite and the edge is used, then
\(I_{p\to q}\ne\emptyset\), so \(N\ne F_p\). If \(N=\emptyset\) then
\(I_{p\to q}=V(F_p)\), which has at least three vertices since \(F_p\)
contains a triangle. Otherwise the boundary edge above gives two
distinct interface vertices \(M',M''\). Either way
\(|I_{p\to q}|\ge 2\). \(\square\)

Thus on a used edge whose source is constrained (bipartite), the
source-side interface is the whole fiber --- both exit colors are
available and there is no exit pin --- and on a used edge at any
non-bipartite fiber at least two interface vertices are available. The
gluing needs slightly more than this bare count when a constrained fiber
exits at a \emph{forced} color, so we must know that color reaches
sufficiently many vertices of the adjacent non-bipartite fiber.

\begin{paperlemma}

\textbf{Lemma 5.7 (delivery avoiding a pinned target).} \emph{Let a
fiber \(C\) be joined to a non-bipartite fiber \(E\) by a used quotient
edge of the buffer line \(L\), the pattern of \(L\) moving a \(1\) from
column \(i\) (pattern \(p\), the \(C\) side) to column \(j\) (pattern
\(q\), the \(E\) side); fix a vertex \(\mathit{bad}\in E\). \textbf{(i)}
If \(C\) is constrained, each color class of \(C\) reaches at least two
distinct vertices of \(E\) along crossing edges; in particular the color
in which \(C\) is forced to exit reaches a vertex \(\ne\mathit{bad}\).
\textbf{(ii)} If \(C\) is non-bipartite, then for every \(x\in C\) some
crossing edge \(yz\) (\(y\in C\), \(z\in E\)) has \(y\ne x\) and
\(z\ne\mathit{bad}\).}

\end{paperlemma}

\emph{Proof.} For \(X\in C\) let \(n_{10}(X),n_{01}(X),n_{11}(X)\) count
the non-pivot rows of type \((1,0),(0,1),(1,1)\) in columns \((i,j)\).
Row \(L\) of \(X\) has type \((1,0)\), contributing \(1\) to column
\(i\) and \(0\) to column \(j\), so counting the \(1\)'s in columns
\(i,j\) over all rows gives \(1+n_{10}(X)+n_{11}(X)=s_i\) and
\(n_{01}(X)+n_{11}(X)=s_j\), whence \[
n_{01}(X)-n_{10}(X)=D:=s_j-s_i+1\qquad(X\in C).\tag{$\ast$}
\] It is worth holding the two sides side by side: a vertex of \(C\) has
pivot row of type \((1,0)\) and imbalance \(D\), a vertex of \(E\) has
pivot row of type \((0,1)\) and imbalance \(2-D\), and in each case the
crossings out of a vertex are counted by its rows of the opposite type.
Each type-\((0,1)\) row of \(X\) yields one crossing from \(X\) into
\(E\) (switch it against row \(L\) on columns \(\{i,j\}\)), distinct
rows giving distinct neighbors, so \(n_{01}(X)\) is the number of
crossings from \(X\) into \(E\). Symmetrically each \(Y\in E\) (whose
row \(L\) has type \((0,1)\)) satisfies \(r(Y)-r'(Y)=s_i-s_j+1=2-D\),
where \(r(Y),r'(Y)\) count its \((1,0)\)- and \((0,1)\)-rows and
\(r(Y)\) is the number of crossings from \(Y\) back to \(C\).

\textbf{(i) \(C\) constrained.} Two values of \(D\) are excluded.

\emph{The case \(D=1\) is impossible.} Then \(s_i=s_j\), so transposing
columns \(i,j\) is an automorphism of \(G(R,S)\) carrying \(p\) to \(q\)
and hence \(C\) isomorphically onto \(E\) --- impossible, since \(C\) is
bipartite and \(E\) is not.

\emph{The case \(D<0\) is impossible.} Then \(2-D\ge 3\), so every
\(Y\in E\) has \(r(Y)\ge 3\); switching \(L\) against three of \(Y\)'s
\((1,0)\)-rows produces three vertices of \(C\) pairwise joined by
interchanges on columns \((i,j)\) --- a triangle in the bipartite \(C\),
again impossible.

Hence \(D=0\) or \(D\ge 2\).

\emph{The case \(D\ge 2\).} Then \((\ast)\) gives \(n_{01}(X)\ge 2\) for
every \(X\in C\), so every vertex of \(C\) carries at least two
crossings into \(E\).

\emph{The case \(D=0\).} Every \(Y\in E\) has \(r(Y)\ge 2\), and \(C\)
(constrained, hence balanced bipartite) has \(|V(C)|\) even. When
\(|V(C)|=2\), switching \(L\) against two \((1,0)\)-rows of any
\(Y\in E\) gives the two vertices of \(C\), so each vertex of \(C\) is
adjacent to all of \(E\), and \(|V(E)|\ge 3\) because \(E\) contains a
triangle. When \(|V(C)|\ge 4\), suppose some color class \(\kappa\)
reached a single \(e\in E\). Every \(M\in\kappa\) has whole-fiber
interface (Lemma 5.6) and therefore crosses only to \(e\), and that
crossing is unique, since distinct type-\((0,1)\) rows of \(M\) would
give distinct neighbors in \(E\); so each such \(M\) is obtained from
\(e\) by undoing that single crossing switch, on rows \(\{L,r_M\}\) and
columns \(\{i,j\}\). Two members \(M,M'\) of \(\kappa\), with witness
rows \(r_M\ne r_{M'}\), are then read off \(e\) in columns \(i,j\):
\[\begin{array}{c|cc|cc|cc}
 & \multicolumn{2}{c|}{e} & \multicolumn{2}{c|}{M} & \multicolumn{2}{c}{M'}\\
\text{row} & i & j & i & j & i & j\\\hline
L & 0 & 1 & 1 & 0 & 1 & 0\\
r_M & 1 & 0 & 0 & 1 & 1 & 0\\
r_{M'} & 1 & 0 & 1 & 0 & 0 & 1
\end{array}\] Thus \(M\) and \(M'\) agree in row \(L\) and differ
exactly on rows \(\{r_M,r_{M'}\}\) in columns \(\{i,j\}\), a single
interchange; hence \(\kappa\) is a clique in the bipartite \(C\) ---
impossible.

In every case each color class reaches at least two vertices of \(E\),
so the forced exit color reaches one \(\ne\mathit{bad}\).

\textbf{(ii) \(C\) non-bipartite.} Put \(c_C:=D\) and \(c_E:=2-D\), so
\(c_C+c_E=2\).

\emph{A symmetry reduction.} The conclusion asks only for a crossing
edge whose \(C\)-endpoint avoids the given vertex \(x\) and whose
\(E\)-endpoint avoids the given vertex \(\mathit{bad}\). The two avoided
vertices play interchangeable roles, and both \(C\) and \(E\) are
non-bipartite here, so the statement is symmetric under
\((C,x)\leftrightarrow(E,\mathit{bad})\), which swaps
\(c_C\leftrightarrow c_E\). Since \(c_C+c_E=2\) we may therefore assume
\(c_C\ge 1\). Both fibers contain a triangle, so \(|C|,|E|\ge 3\).

\emph{The case \(c_C\ge 2\).} Then \((\ast)\) gives \(n_{01}(X)\ge 2\)
for every \(X\in C\). Choose \(X\ne x\), possible as \(|C|\ge 3\); then
\(X\) has two crossings to distinct vertices of \(E\), at most one of
them \(\mathit{bad}\), and we take \(y:=X\) and \(z\ne\mathit{bad}\).

\emph{The case \(c_C=1\)}, so that \(c_E=1\) and \(s_i=s_j\). Suppose
every crossing were incident to \(x\) or to \(\mathit{bad}\). Pick
distinct \(X_1,X_2\in C\setminus\{x\}\); each has a crossing, since
\(n_{01}(X_t)=n_{10}(X_t)+1\ge 1\), and that crossing does not run to
\(x\), hence runs to \(\mathit{bad}\), so \(r(\mathit{bad})\ge 2\). The
crossing \(\mathit{bad}\,X_1\) switches rows \(\{L,\rho\}\) for a
\((1,0)\)-row \(\rho\) of \(\mathit{bad}\), so the type-\((0,1)\) rows
of \(X_1\) are those of \(\mathit{bad}\) together with \(\rho\), giving
\(n_{01}(X_1)=r'(\mathit{bad})+1=r(\mathit{bad})\ge 2\), using
\(r(\mathit{bad})-r'(\mathit{bad})=2-D=1\). But then every crossing of
\(X_1\) runs to \(\mathit{bad}\), so \(n_{01}(X_1)\le 1\) --- a
contradiction.

In either case a crossing \(yz\) has \(y\ne x\) and
\(z\ne\mathit{bad}\); a Hamilton path of \(C\) from \(x\) to \(y\)
(Hamilton-connectedness, induction hypothesis) then delivers
\(z\ne\mathit{bad}\) into \(E\). \(\square\)

Clause (i) lets a constrained fiber, whose exit color is forced, still
hand the non-bipartite target an entry distinct from its pinned exit;
clause (ii) is the symmetric non-bipartite-source case. A singleton
source is a deterministic pass-through: its used edge is already
realized, and the two interface vertices of Lemma 5.6 on the
non-bipartite side supply an entry \(\ne\mathit{bad}\).

\hypertarget{a-separating-buffer-line-lemmas-5.85.10}{%
\subsubsection{A Separating Buffer Line (Lemmas
5.8--5.10)}\label{a-separating-buffer-line-lemmas-5.85.10}}

\hypertarget{triangles-and-wide-pairs}{%
\paragraph{Triangles and wide pairs}\label{triangles-and-wide-pairs}}

Call a pair of rows a \textbf{wide pair} if, in some realization of
\(G(R,S)\), their supports are incomparable (each has a \(1\) in some
column where the other has a \(0\)) and their symmetric difference has
size at least \(3\); \textbf{wide column pairs} are defined dually.
Triangles are governed by these.

\begin{paperlemma}

\textbf{Lemma 5.8 (triangle classification).} \emph{Every triangle of
\(G(R,S)\) is supported on exactly two rows and three columns --- its
two rows then forming a wide pair --- or, dually, on exactly three rows
and two columns, its two columns forming a wide column pair.}

\end{paperlemma}

\emph{Proof.} Throughout this proof \(M\triangle M'\) denotes the set of
cells at which the two matrices differ; since the matrices are
\(0\)-\(1\), this is the symmetric difference of their supports, and two
matrices are adjacent in \(G(R,S)\) exactly when \(M\triangle M'\) is
the four-cell support of a single interchange. Two such sets combine as
sets do:
\(M_0\triangle M_2=(M_0\triangle M_1)\triangle(M_1\triangle M_2)\),
because a cell differs between \(M_0\) and \(M_2\) precisely when it
differs an odd number of times along the way.

Let \(M_0,M_1,M_2\) be a triangle. Each symmetric difference
\(M_i\triangle M_j\) is the support of a single interchange: four cells
\(\{a,b\}\times\{i,j\}\) carrying the \(2\times2\) checkerboard pattern.
From \(M_0\triangle M_2=(M_0\triangle M_1)\triangle(M_1\triangle M_2)\),
the symmetric difference of the two interchange-supports
\(S_1=M_0\triangle M_1\) and \(S_2=M_1\triangle M_2\) is itself a single
interchange-support, so \(|S_1\triangle S_2|=4\) and hence
\(|S_1\cap S_2|=2\). Writing \(S_1=\{a,b\}\times\{i,j\}\), the two
shared cells lie in
\((\{a,b\}\cap\mathrm{rows}(S_2))\times(\{i,j\}\cap\mathrm{cols}(S_2))\);
two cells force either both rows shared with a single column, or a
single row shared with both columns (a shared diagonal would force
\(S_2=S_1\)). The exact complementary cardinality follows directly. In
the first case the two supports use the same two rows, and their
two-column sets intersect in one column; since \(S_1\ne S_2\), their
union therefore has \(2+2-1=3\) columns. In the second case their
two-row sets intersect in one row and their column sets agree, so their
union has \(2+2-1=3\) rows and exactly two columns. The remaining
support \(M_0\triangle M_2\) lies in the same two rows (resp. columns),
so the whole triangle is supported on two rows and three columns, or on
three rows and two columns. Its two rows (resp. columns) carry the three
rotations of the \(2\times3\) (resp. \(3\times2\)) pattern, so they
differ in all three changing columns (rows) --- a wide pair. \(\square\)

Figure 3 shows the two types.

\begin{figure}
\centering
\includegraphics[width=0.88\textwidth,height=\textheight]{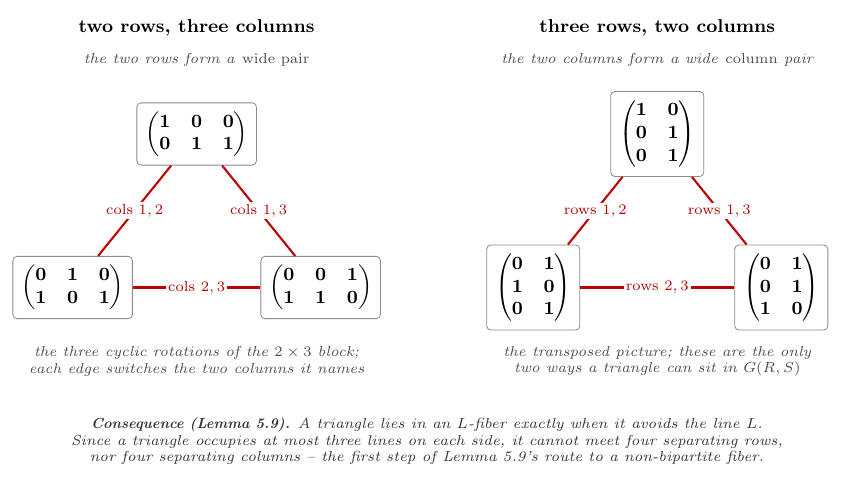}
\caption{\textbf{The two triangle types of Lemma 5.8.} A triangle of
\(G(R,S)\) occupies exactly two rows and three columns, or exactly three
rows and two columns. Each is drawn as the three cyclic rotations of a
\(2\times3\) (respectively \(3\times2\)) block, every edge labeled by
the pair of columns (rows) whose interchange it performs. In both
triples the row and column sums agree across the three matrices and each
pair differs by a single \(2\times2\) interchange. A triangle lies in an
\(L\)-fiber exactly when it avoids the line \(L\); since a triangle
occupies at most three lines on each side, it cannot meet four
separating rows, nor four separating columns. That is the first step of
Lemma 5.9's argument that some separating line has a non-bipartite
fiber, and the remaining cases need the wide-pair analysis given there.}
\end{figure}

Since a non-bipartite interchange graph contains a triangle (Lemma 2.1,
equivalently Brualdi \citep[Thm.
6.3.4]{brualdiCombinatorialMatrixClasses2006}), it has a wide row pair
or a wide column pair.

\hypertarget{buffer-existence}{%
\paragraph{Buffer existence}\label{buffer-existence}}

\begin{paperlemma}

\textbf{Lemma 5.9 (buffer existence).}\\
\emph{For every pair of distinct vertices \(a,b\) in the present
non-base case, there exists a line \(L\) separating \(a\) from \(b\)
such that one fiber of \(L\) is non-bipartite.}

\end{paperlemma}

\emph{Proof.} A triangle of \(G\) lies in an \(L\)-fiber exactly when it
does not use the line \(L\). Each fiber is itself an interchange graph,
on the class obtained by deleting \(L\) (§2.1), so Lemma 2.1 says that a
fiber is non-bipartite exactly when it contains a triangle. Hence a
separating line \(L\) has a non-bipartite fiber if and only if some
triangle avoids \(L\).

Suppose no separating line has a non-bipartite fiber. Then every
triangle uses every separating row and every separating column of the
pair \(a,b\). Let \(R^*\) and \(C^*\) be the separating rows and
columns. Since \(a\) and \(b\) have the same margins,
\(|R^*|,|C^*|\ge2\). Indeed, form the bipartite graph on the row and
column indices whose edges are the cells where \(a\) and \(b\) differ,
coloring an edge according to whether the entry is \(1\) in \(a\) or in
\(b\). Equality of the margins gives the two colors equal degree at
every vertex, so these edges decompose into color-alternating cycles;
each such cycle meets at least two rows and at least two columns, and
every line it meets is separating.

Every triangle in an interchange graph is either a two-row triangle
using exactly two rows and three columns, or a three-row triangle using
exactly three rows and two columns (Lemma 5.8). Therefore, if
\(|R^*|\ge4\) or \(|C^*|\ge4\), no triangle can meet all separating
lines, contradiction. If \(|R^*|=|C^*|=3\), a two-row triangle misses a
separating row and a three-row triangle misses a separating column,
again a contradiction. In summary:

\[\begin{array}{ll}
|R^*|\ge4\text{ or }|C^*|\ge4 & \text{impossible: every triangle misses a separating line}\\
(3,3) & \text{impossible: each triangle type misses one side}\\
(2,2),\ (2,3),\ (3,2) & \text{the remaining cases, handled below}
\end{array}\]

We also use the converse of the triangle classification: a wide row pair
\(\{p,q\}\) supports a triangle on those two rows. Their symmetric
difference \(\Delta\) has \(|\Delta|\ge3\) with both sides nonempty, so
one side---say \(p\)'s---has at least two columns; taking two columns
from \(p\)'s side and one from \(q\)'s gives a \(2\times3\) block whose
three rotations are interchanges (each preserves all margins), hence lie
in \(G\) and form a triangle on rows \(\{p,q\}\) (dually for columns).

In the remaining cases, every wide row pair must be the unique two-row
set \(R^*\) when \(|R^*|=2\), and no wide row pair can exist when
\(|R^*|=3\). Similarly, every wide column pair must be the unique
two-column set \(C^*\) when \(|C^*|=2\), and no wide column pair can
exist when \(|C^*|=3\). Thus \(G\) has at most one wide row pair and at
most one wide column pair.

\textbf{Claim: a unique wide row pair (or column pair) cannot exist.}
Within this claim \(p,q,x,z\) denote row indices --- not patterns, as
elsewhere in this section --- and \(R_y(N)\) denotes the support of row
\(y\) in a realization \(N\). Suppose \(\{p,q\}\) is the unique wide row
pair, and let \(\mathcal C\) be the unique wide column pair if one
exists, and \(\mathcal C=\emptyset\) otherwise. Choose a realization
\(M\) in which \(p,q\) are wide, and write \[
\Delta_p=R_p(M)\setminus R_q(M),\quad \Delta_q=R_q(M)\setminus R_p(M),\quad \Delta=\Delta_p\cup\Delta_q,
\] with \(\Delta_p,\Delta_q\ne\emptyset\) and \(|\Delta|\ge3\). Choose
\(d\in \Delta\) not belonging to \(\mathcal C\), possible because
\(\mathcal C\) has at most two columns.

\emph{Step 1: no switch can use exactly one of \(p,q\).} Along any
interchange path from \(M\), track the current supports of rows \(p,q\)
through the split \[
\Delta_p^t=R_p\setminus R_q,\qquad \Delta_q^t=R_q\setminus R_p,\qquad (\Delta_p^0,\Delta_q^0)=(\Delta_p,\Delta_q),
\] where \(R_p,R_q\) denote the supports in the current realization, and
consider the invariant ``\emph{rows \(p,q\) are wide with the same
difference set \(\Delta=\Delta_p^t\cup\Delta_q^t\)}'': the split
\((\Delta_p^t,\Delta_q^t)\) may change from step to step, but \(\Delta\)
and the side sizes do not.

A switch avoiding both \(p\) and \(q\) leaves their supports unchanged,
and a switch using both exchanges one column of \(\Delta_p^t\) with one
of \(\Delta_q^t\), preserving \(\Delta\) and the two side sizes; so the
invariant can first fail only at a switch using exactly one of \(p,q\).
Take the \emph{first} such switch (if any) --- one that uses exactly one
row from \(\{p,q\}\) and one outside row --- say rows \(p\) and
\(x\notin\{p,q\}\), on columns \(u,v\). Just before it the invariant
still holds. The switch's checkerboard puts
\(u\in\operatorname{supp}(p)\setminus\operatorname{supp}(x)\) and
\(v\in\operatorname{supp}(x)\setminus\operatorname{supp}(p)\), so
\(p,x\) are incomparable in the current realization; since \(\{p,x\}\)
is not wide, their symmetric difference is exactly \(\{u,v\}\). The two
rows therefore agree outside \(u,v\), and the switch exchanges them
there: \[\begin{array}{c|cc}
\text{row} & u & v\\\hline
p\ \text{before} & 1 & 0\\
x\ \text{before} & 0 & 1\\\hline
p\ \text{after} & 0 & 1\\
x\ \text{after} & 1 & 0
\end{array}\] So after the switch \(R_x\) equals row \(p\)'s pre-switch
support: \(\{x,q\}\) inherits the pre-switch split
\((\Delta_p^t,\Delta_q^t)\) of \(\Delta\) and is wide in the post-switch
realization --- a second wide row pair, contradiction. Thus along any
path \(p,q\) can only switch with each other, and by the previous
paragraph the invariant persists.

\emph{Step 2: pivot-freedom creates a second wide column pair.} We have
shown that every path out of \(M\) preserves \(\Delta\) and that no
switch on such a path uses exactly one of \(p,q\). We now force a change
in a cell outside those two rows, and show that it creates a second wide
column pair.

Pick a third row \(x\notin\{p,q\}\). By Lemma 5.2, the cell \((x,d)\)
varies; choose a matrix in which its value differs from that in \(M\).
By Ryser connectivity (§1), join \(M\) to that matrix by an interchange
path and take the first switch on the path that changes \(x_d\). By the
previous paragraph this switch uses only rows outside \(\{p,q\}\). Just
before the switch, choose \(c\in \Delta\) on the opposite side of the
current \(p/q\) split from \(d\). Then columns \(c,d\) differ on rows
\(p,q\) with opposite orientations. Since \(d\notin\mathcal C\), the
column pair \(\{c,d\}\) is not wide. Therefore every outside row \(z\)
has \(z_c=z_d\); otherwise columns \(c,d\) would differ on \(p,q,z\),
forming a wide column pair. In particular \(x_c=x_d\), so the first
switch changing \(x_d\) cannot use column \(c\); writing \(\varepsilon\)
for that common value, and taking \(d\) on the \(p\) side of the split
and \(c\) on the \(q\) side, \[\begin{array}{c|cc}
\text{row} & c & d\\\hline
p & 0 & 1\\
q & 1 & 0\\\hline
x\ \text{before} & \varepsilon & \varepsilon\\
x\ \text{after} & \varepsilon & 1-\varepsilon
\end{array}\] After the switch \(x_c\ne x_d\), so columns \(c,d\) differ
on all three of \(p,q,x\) --- a wide column pair other than
\(\mathcal C\), a contradiction.

Thus no such unique wide row pair exists; by transposition no unique
wide column pair exists. Since \(G\) is non-bipartite, Lemma 5.8 and its
converse give at least one wide row or column pair. Contradiction. Hence
some separating line has a non-bipartite fiber. \(\square\)

\begin{paperlemma}

\textbf{Lemma 5.10 (every line has non-bipartite quotient).} \emph{In a
nonempty active invariant-free class with \(m,n\ge3\), the quotient
\(Q_L\) is non-bipartite for every line \(L\). In particular, this holds
for the buffer line supplied by Lemma 5.9.}

\end{paperlemma}

\emph{Proof.} Let \(L\) be a row; the column case follows by
transposition. The realizable patterns of \(L\) form a \textbf{shifted}
matroid \(M_L\): order the columns so \(s_1\ge s_2\ge\cdots\); writing
\(\sigma(X)=\sum_{j\in X}s_j\) and \(p(X)=|p\cap X|\), if \(v\in p\),
\(u\notin p\) and \(s_u\ge s_v\), then for every \(X\) with
\(v\in X,\,u\notin X\) the set \(X'=X-v+u\) satisfies
\(\sigma(X')\ge\sigma(X)\) --- hence \(\ell(X')\ge\ell(X)\), since \(h\)
depends only on the cardinality --- and \(p'(X)=p(X')\), so \(p'=p-v+u\)
still meets the Gale--Ryser inequalities and is realizable. (For a set
\(X\) with \(u\in X\) or \(v\notin X\), the shift leaves \(|p\cap X|\)
unchanged or larger, so that inequality is preserved automatically; only
\(u\notin X,\,v\in X\) needs the argument just given.) In an active
invariant-free block every cell of \(L\) is variable (Lemma 5.2), so
\(M_L\) is loopless and coloopless on its \(f\ge3\) active columns
(\(f=n\) for a row line, \(f=m\) for a column line; \(f\ge3\) by the
regime), of some rank \(k\) with \(1\le k\le f-1\).

We use the following elementary consequence of shiftedness, and call the
process \textbf{compression}. If a basis is required to contain a fixed
element \(e\), or required to avoid \(e\), then repeated downward shifts
that preserve that requirement carry it to the lexicographically first
\(k\)-set satisfying the same requirement, and every intermediate set is
a basis. Indeed, let \(T\) be that first \(k\)-set. If a current basis
\(B\ne T\) satisfies the requirement, choose the least
\(u\in T\setminus B\). Some \(v\in B\setminus T\) has \(v>u\): otherwise
every element of \(B\setminus T\) would precede \(u\), making the least
element of \(B\triangle T\) lie in \(B\), so that \(B\) itself would be
lexicographically earlier than \(T\) --- contradicting \(T\)'s choice as
the lexicographically first \(k\)-set satisfying the requirement.
Replacing \(v\) by \(u\) is a permitted downward shift and preserves the
requirement. Iteration terminates at \(T\).

If \(k=1\), looplessness gives the singleton bases
\(\{1\},\{2\},\{3\}\), which form a triangle in the base graph. Assume
\(k\ge2\). With no constraint, compression gives \(B_0=\{1,\dots,k\}\).
Since \(k+1\) is not a loop, some basis contains \(k+1\); compressing
among bases containing \(k+1\) gives \(B_1=\{1,\dots,k-1,k+1\}\). Since
\(k-1\) is not a coloop, some basis omits \(k-1\); compressing among
bases omitting \(k-1\) gives \(B_2=\{1,\dots,k-2,k,k+1\}\). The three
bases \(B_0,B_1,B_2\) are pairwise related by a single exchange, so the
base graph contains a triangle. By Lemma 5.4, \(Q_L\) is the base graph,
so it contains this triangle and is non-bipartite. (Here \(U_{r,f}\)
denotes the uniform matroid of rank \(r\) on \(f\) elements and
\(\oplus\) the direct sum of matroids. The Boolean square
\(U_{1,2}\oplus U_{1,2}\), whose base graph is the bipartite
\(4\)-cycle, is \emph{not} shifted, so it never arises as an \(M_L\).)
\(\square\)

\hypertarget{gluing-the-fibers-proposition-5.11}{%
\subsubsection{Gluing the fibers (Proposition
5.11)}\label{gluing-the-fibers-proposition-5.11}}

Figure 4 shows the interior-buffer case of the construction.

\begin{paperproposition}

\textbf{Proposition 5.11 (single-pass gluing).} \emph{Let \(L\) separate
distinct vertices \(a,b\). Assume, as supplied here by the induction
hypothesis of §3, that every non-singleton fiber of \(L\) is
Hamilton-laceable when bipartite and Hamilton-connected when
non-bipartite. If (1) \(Q_L\) is non-bipartite, hence Hamilton-connected
by Corollary 5.5, and (2) some fiber of \(L\) is non-bipartite, then
\(G\) has a Hamilton path from \(a\) to \(b\) visiting each fiber once.}

\end{paperproposition}

\begin{figure}
\centering
\includegraphics[width=0.95\textwidth,height=\textheight]{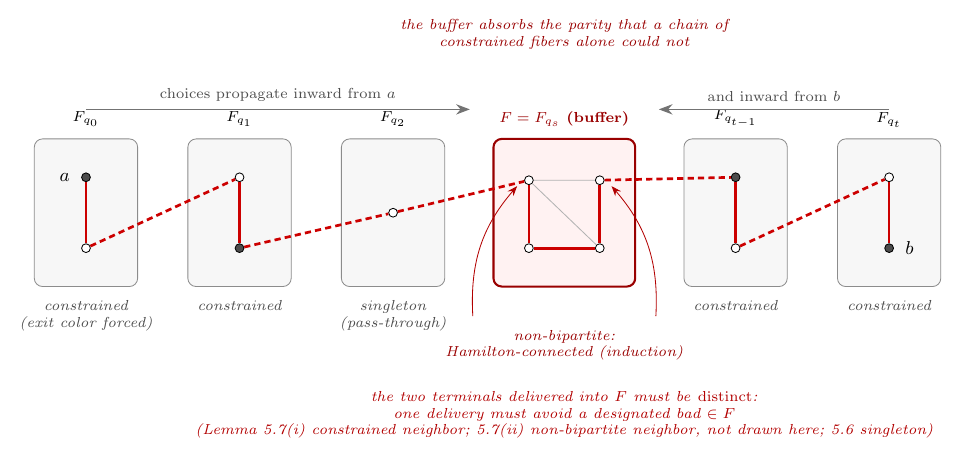}
\caption{\textbf{Single-pass gluing, in the interior-buffer case.} The
branch Figure 2 does not show, the worked example closing this section
being the endpoint-buffer case; the drawing is schematic rather than an
instance. The fibers appear in the order of the chosen Hamilton path of
\(Q_L\), and only that path, the chosen internal traversals and the
chosen crossing edges are drawn --- all other edges of \(G\), and all
quotient edges off the path, are suppressed. Solid red edges are
internal traversals and dashed red edges are splices; together they form
the Hamilton \(a\)--\(b\) path. Inside a constrained fiber, black and
white mark its two local color classes, the names not being coordinated
between fibers. The flanking fibers drawn are representative: a
constrained fiber exits in the forced color, a singleton is a
deterministic pass-through, and a non-buffer non-bipartite fiber may
also occur, though none is drawn here. Terminal choices propagate inward
from the pinned endpoints \(a\) and \(b\), so the two entries delivered
into the buffer \(F\) must be distinct; which lemma supplies that
distinctness depends on the delivering fiber, and the proof below
settles it case by case.}
\end{figure}

\emph{Proof.} Because \(L\) is separating, each of its fibers is a
proper subclass and hence, by §2.1, an interchange graph with strictly
fewer vertices than \(G\); after the canonical normalization of Lemma
5.1 the induction hypothesis therefore applies to it. With the bipartite
branch supplied by Lemma 1.2, this gives a Hamilton-laceable traversal
of every non-singleton bipartite fiber and a Hamilton-connected
traversal of every non-bipartite one. Let \(\alpha,\beta\) be the
quotient vertices holding \(a,b\), and choose a Hamilton path
\(\alpha=q_0,\ldots,q_t=\beta\) of \(Q_L\) (Hamilton-connected). Since
\(Q_L\) is non-bipartite this path has at least three vertices. The
non-bipartite fiber \(F=F_{q_s}\) of hypothesis (2) lies on the path. We
traverse each fiber by one internal path and splice consecutive fibers
on crossing edges.

Orient each flanking chain from its pinned endpoint toward \(F\),
reversing the displayed order \(q_0,\ldots,q_t\) on the far side where
necessary. Propagation then maintains one invariant: the current fiber
has a prescribed entry, we choose an admissible exit, and a crossing
edge out of that exit fixes the entry of the next fiber. Every interface
lemma below is applied in this direction.

\emph{Fiber traversals.} A constrained fiber, once its entry is fixed,
exits in the forced opposite color; its whole-fiber interface (Lemma
5.6) realizes an outgoing crossing from any exit vertex of that color,
so the next entry may be pinned but this is no obstruction --- that
fiber simply exits opposite. A singleton is a deterministic
pass-through. A non-bipartite fiber other than \(F\) is traversed
one-directionally: its entry is delivered from upstream and its exit is
chosen on the interface toward the next fiber, which offers at least two
vertices (Lemma 5.6), one distinct from the entry;
Hamilton-connectedness (induction hypothesis) joins them.

\emph{The buffer.} Remove \(F\) from the quotient path; the remaining
vertices form the chains flanking \(q_s\) --- two chains when \(F\) is
interior (\(0<s<t\)), one when \(F\) is an endpoint; we call these the
chains \emph{flanking} \(F\). Propagate terminal choices inward from the
pinned endpoints \(a,b\) along each chain toward \(F\), in the following
order. If \(F\) is interior, propagate one flanking chain first and let
\(\mathit{bad}\) be the entry into \(F\) that it delivers; then
propagate the other chain, whose delivered entry must differ from
\(\mathit{bad}\). If \(F\) is an endpoint, take \(\mathit{bad}\) to be
the prescribed endpoint of \(G\) lying in \(F\), and propagate the
single flanking chain. The two terminals of \(F\) to be joined are then:
when \(F\) is interior, the two entries delivered by the flanking
chains, which must be distinct; when \(F\) is an endpoint, the entry
delivered by its one flanking chain together with the pinned endpoint of
\(G\) lying in \(F\), which the entry must avoid. Either way one
delivery into \(F\) must avoid a designated vertex
\(\mathit{bad}\in F\), and Lemma 5.7 supplies it --- clause (i) when the
delivering neighbor is constrained, clause (ii) when it is
non-bipartite, and Lemma 5.6 when it is a singleton. With the two
distinct terminals fixed, Hamilton-connectedness of \(F\) (induction
hypothesis; \(|V(F)|\ge 3\)) joins them through all of \(F\).

Splicing every internal path along the selected crossing edges yields a
single Hamilton path from \(a\) to \(b\) using each vertex of each fiber
exactly once. The non-bipartite buffer absorbs the parity that a chain
of constrained fibers alone could not. \(\square\)

Combining Lemma 5.9 (a non-bipartite buffer fiber) with Lemma 5.10 (the
same line has non-bipartite, hence Hamilton-connected, quotient \(Q_L\))
yields a separating line satisfying both hypotheses of Proposition 5.11.
Therefore \(G\) has a Hamilton path from \(a\) to \(b\), completing the
indecomposable non-base case.

\hypertarget{the-construction-as-an-algorithm}{%
\paragraph{The construction as an
algorithm}\label{the-construction-as-an-algorithm}}

For completeness, and to mirror the computational verification, we state
the single-pass construction explicitly; each step is licensed by the
lemma named beside it.

\emph{Input.} An active, invariant-free, non-bipartite, indecomposable
class \(G=G(R,S)\) with \(m,n\ge 3\), and distinct endpoints \(a,b\).

\begin{enumerate}
\def\labelenumi{\arabic{enumi}.}
\tightlist
\item
  \textbf{(Lemma 5.9)} Choose a separating line \(L\) --- one with
  \(L(a)\ne L(b)\) --- having a non-bipartite fiber (the \emph{buffer
  line}); such a line exists.
\item
  \textbf{(Lemmas 5.3, 5.4, 5.10; Corollary 5.5)} Its quotient \(Q_L\)
  is non-bipartite (5.10); the patterns of \(L\) are the bases of a
  matroid (5.3) and \(Q_L\) is that matroid's base-exchange graph (5.4),
  so \(Q_L\) is a non-bipartite matroid base-exchange graph and hence
  Hamilton-connected (Naddef--Pulleyblank, Corollary 5.5).
\item
  Take a Hamilton path \(\alpha=q_0,q_1,\dots,q_t=\beta\) of \(Q_L\)
  between the patterns \(\alpha,\beta\) of \(a,b\) (it exists by step
  2).
\item
  \textbf{(Lemmas 5.6, 5.7; Proposition 5.11)} Thread the fibers along
  \(q_0,\dots,q_t\): traverse each fiber by a Hamilton path obtained by
  recursion on that strictly smaller interchange class, crossing between
  consecutive fibers on a selected interface edge. At each ordinary
  splice incident to the non-bipartite buffer fiber there are at least
  two interface vertices (Lemma 5.6); at the splice that must avoid a
  vertex pinned from the other direction, the required freedom is
  supplied by Lemma 5.7 --- clause (i) for a constrained delivering
  fiber, clause (ii) for a non-bipartite one --- or by Lemma 5.6 for a
  singleton. The interior- and endpoint-buffer cases are the two
  branches of Proposition 5.11.
\end{enumerate}

\emph{Output.} A Hamilton \(a\)--\(b\) path of \(G\) (Proposition 5.11).

\hypertarget{verification}{%
\paragraph{Verification}\label{verification}}

The construction above is exactly the procedure executed by the verifier
\texttt{verify\_capstone\_faithful.py}, which drives each step by the
lemma it cites and aborts if any cited guarantee fails. A companion
verifier, \texttt{verify\_capstone\_deterministic.py}, runs the same
construction with \emph{every} interface choice fixed by the rule the
proof prescribes, with no search. It produced a valid Hamilton path on
the deterministic first pass, with \textbf{zero backtracks and zero
repairs}, for every endpoint pair of every §5 pivot class within the
computational bound; both the interior-buffer and endpoint-buffer
branches of the gluing were exercised. The construction is therefore
genuinely deterministic, not merely existence-certified. On every
non-bipartite class within the computational bound (exhaustively over
endpoint pairs for the smaller classes, and a sampled sweep of the
larger ones), every cited guarantee held and the procedure returned a
valid Hamilton path for every pair tested. The verifiers and their
ledger are released with the paper's verification artifacts.

\hypertarget{a-worked-example}{%
\paragraph{A worked example}\label{a-worked-example}}

We close the section by running the construction once, in full, on a
small class with no singleton fibers and with three constrained fibers,
so every quotient step before the buffer is forced: \(R=(2,2,1)\),
\(S=(2,1,1,1)\). This is the class drawn in Figure 2, which the reader
may wish to have alongside what follows: the twelve matrices, the four
fibers, the quotient \(K_4\) and the lifted path all appear there. It is
not the smallest mixed-fiber class; for example, \(R=(3,2,1)\),
\(S=(2,2,1,1)\) has eight matrices and fiber sizes \(1,2,5\) for a
column of sum two.

Row \(3\) has a single \(1\), in some column \(c\). If \(c\ne1\), then
column \(1\) (of sum \(2\)) lies in both rows \(1\) and \(2\), and the
two remaining sum-one columns split between them; if \(c=1\), the
residual column sums are \((1,1,1,1)\) and rows \(1,2\) take
complementary \(2\)-subsets of the columns. Thus a matrix is determined
by the pair \((c\mid xy)\), where \(xy\) is the support of row \(1\);
for instance \[
a=\begin{pmatrix}1&1&0&0\\1&0&1&0\\0&0&0&1\end{pmatrix}=(4\mid 12),
\qquad
b=\begin{pmatrix}1&1&0&0\\0&0&1&1\\1&0&0&0\end{pmatrix}=(1\mid 12).
\] There are \(12\) matrices; every cell takes both values across
\(\mathcal A(R,S)\), so the class is active and invariant-free, hence
prime (§2.2), non-bipartite (the fiber \(F_1\) below contains
triangles), with \(m,n\ge3\) and \(|V|>6\): the §5 regime. Fix the
endpoints \(a,b\) displayed above.

\emph{The buffer line (Lemma 5.9).} Take \(L=\) row \(3\); it separates
\(a\) from \(b\) (their patterns are \(\{4\}\) and \(\{1\}\)). Its
fibers, indexed by the pattern: \(F_1\) consists of the six matrices
\((1\mid xy)\), and a single interchange moves one element of the
row-\(1\) support, so \(F_1\) is the Johnson graph \(J(4,2)\) --- the
octahedron, non-bipartite: the \textbf{buffer}. \(F_2,F_3,F_4\) have two
matrices each (the two splits of the free columns), one interchange
apart: bipartite \(K_2\)'s --- \textbf{constrained} fibers.

\emph{The quotient (Lemmas 5.3, 5.4, 5.10; Corollary 5.5).} The
realizable patterns are all four singletons --- the bases of the uniform
matroid \(U_{1,4}\) --- and every exchange is realized by a crossing
edge (Lemma 5.4: moving row \(3\)'s \(1\) between two columns needs a
witness row of type \((0,1)\) there, and one is always present). So
\(Q_L\) is the base-exchange graph \(K_4\): non-bipartite and
Hamilton-connected --- here visibly, in general by Corollary 5.5.

\emph{The lift (Proposition 5.11, the endpoint-buffer case).} Choose the
Hamilton path \(\{4\},\{3\},\{2\},\{1\}\) of \(Q_L\) from \(a\)'s
pattern to \(b\)'s; the non-bipartite fiber sits at the \(b\) end.
Thread: \[\begin{aligned}
&\underbrace{(4\mid 12)\,(4\mid 13)}_{F_4}\;\big|\;
\underbrace{(3\mid 14)\,(3\mid 12)}_{F_3}\;\big|\;
\underbrace{(2\mid 13)\,(2\mid 14)}_{F_2}\;\big|\;\\
&\underbrace{(1\mid 24)\,(1\mid 34)\,(1\mid 14)\,(1\mid 13)\,(1\mid 23)\,(1\mid 12)}_{F_1}
\end{aligned}\] Each \(K_2\) is traversed by its edge, so its entry
vertex forces its exit vertex; the crossing out of the forced exit
nonetheless always exists, because a constrained fiber on a used edge
has whole-fiber interface (Lemma 5.6 --- indeed both vertices of each
\(K_2\) here carry crossings into the next fiber). Each crossing is one
interchange moving row \(3\)'s \(1\) against a witness row: for example
\((4\mid 13)\to(3\mid 14)\) switches rows \(\{1,3\}\) on columns
\(\{3,4\}\).

The one delicate step is the final crossing, into the buffer, which must
avoid the pinned endpoint \(b\). Here the counting identity of Lemma 5.7
gives \(D=s_1-s_2+1=2\), so \emph{every} vertex of \(F_2\) has at least
two crossings into \(F_1\) --- the \(D\ge2\) branch: from the forced
exit \((2\mid 14)\) the two targets are \((1\mid 24)\) and
\((1\mid 14)\), and we take \((1\mid 24)\ne b\). Inside the buffer,
\(F_1\) is a strictly smaller interchange class, non-bipartite, hence
Hamilton-connected by the induction hypothesis (concretely: the
octahedron, a §6 base object), and the displayed path runs from the
entry \((1\mid 24)\) to \(b=(1\mid 12)\) through all six vertices.

The result is a Hamilton \(a\)--\(b\) path of \(G\) visiting each of the
\(12\) matrices exactly once. The example also shows why the buffer is
necessary: in the three \(K_2\) fibers every internal traversal is
forced, so a chain consisting only of constrained fibers would deliver a
\emph{predetermined} final vertex --- the wrong one for half the
demands. The non-bipartite fiber restores the choice; that is the parity
absorption of Proposition 5.11 in miniature. (Had both endpoints lain in
constrained fibers, the same lift would run with the buffer as an
interior fiber, via Proposition 5.11 (interior-buffer case).)

\begin{center}\rule{0.5\linewidth}{0.5pt}\end{center}

\hypertarget{base-classes}{%
\subsection{6. Base classes}\label{base-classes}}

This section discharges the base cases, case (ii) of the non-bipartite
branch, together with the smallest bipartite blocks, all without the
induction hypothesis: the single complete transposition graphs, the
Johnson graphs, and the classes of order at most six. Each is settled
directly or by a cited theorem.

A single complete transposition graph is Hamilton-laceable by the
Bipartite 2-DPC Lemma (§7.2); independently, for \(a\ge4\), the stronger
theorem of Araki \citep{arakiHyperHamiltonianLaceability2006a} on Cayley
graphs of \(S_a\) with a connected transposition generating set implies
the same Hamilton-laceability conclusion without the Coleman input.

If at most two rows of \(G\) are active, or at most two columns, then
\(G\) is isomorphic to a Johnson graph \(J(f,k)\), which is
Hamilton-connected by Alspach \citep{alspachJohnsonGraphsAre2012}.
Indeed, the class here has already been normalized by the deletion of
inactive lines (§3), so every one of its lines is active; at most two
active rows therefore means that the normalized class has exactly two
rows, a single row leaving no column active. Write \(f\) for its number
of columns and \(r_1\) for the first of its two row sums. Every column
has sum \(s_j\) with \(0<s_j<2\), hence \(s_j=1\), so each column
carries its single \(1\) in row \(1\) or in row \(2\). A matrix of the
class is therefore determined by the set of columns whose \(1\) lies in
row \(1\), an \(r_1\)-subset of the \(f\) columns, and an interchange on
rows \(\{1,2\}\) and columns \(\{j,k\}\) exchanges one element of that
set for one outside it. This identifies \(G\) with \(J(f,r_1)\); the
column case is dual. (Both \(f\) and \(r_1\) are read off the normalized
class: deleting an inactive column, all-zero or all-ones, lowers \(f\)
by one, and an all-ones column also lowers every row sum by one.)

Finally, the interchange graphs of order at most \(6\) are checked
directly. The list of cases is complete because a nonempty prime active
class has no invariant positions by the Brualdi--Manber prime criterion
(§2.2), and Brualdi's assignment-polytope count \citep[§9.13,
pp.~491--493, footnote 26, after
Brualdi--Hartfiel--Hwang]{brualdiCombinatorialMatrixClasses2006}
therefore gives \(|\mathcal A(R,S)|\ge(m-1)(n-1)+1\), so order \(\le 6\)
forces at most three active rows and three active columns beyond the
two-line Johnson cases. An active \(3\times3\) margin has all line sums
in \(\{1,2\}\). Up to permutations this leaves the two bipartite
constant classes \(R=S=(1,1,1)\) and \(R=S=(2,2,2)\), together with
\(R=S=(2,1,1)\) and its complement \(R=S=(2,2,1)\).

Each of the latter two is a five-matrix class whose interchange graph is
\(K_5\) minus two disjoint edges. That graph is Hamilton-connected on
every endpoint pair: label the two missing edges \(x_1x_2\) and
\(y_1y_2\) and let \(z\) be the fifth vertex, so that up to automorphism
the endpoint pairs are represented by \((x_1,x_2)\), \((x_1,y_1)\) and
\((z,x_1)\), with Hamilton paths \(x_1y_1zy_2x_2\), \(x_1y_2x_2zy_1\)
and \(zy_1x_2y_2x_1\) respectively. (The formal companion checks all
pairs by an executable certificate.)

This completes the reduction: assuming the Bipartite 2-DPC Lemma, every
interchange graph is maximally Hamiltonian.

\begin{center}\rule{0.5\linewidth}{0.5pt}\end{center}

\hypertarget{proof-of-the-bipartite-2-dpc-lemma}{%
\subsection{7. Proof of the Bipartite 2-DPC
Lemma}\label{proof-of-the-bipartite-2-dpc-lemma}}

This section proves the Bipartite 2-DPC Lemma (Lemma 1.2), the
paired-cover fact that the bipartite case of Section 3 and the
one-bipartite-factor case of Section 4 both consume. It stands outside
the main induction. The proof here inducts on the number of Cartesian
factors; the recursion on rank is internal to the
transposition-like-graph theorem it uses as an external input.

By Proposition 2.2, the Bipartite 2-DPC Lemma reduces to the following
(orders below \(4\) are vacuous: a paired \(2\)-disjoint-path-cover
demand needs four distinct vertices, so \(K_2=CT_2\) satisfies the
Bipartite 2-DPC Lemma trivially).

\begin{papertheorem}

\textbf{Theorem 7.1.} \emph{Every Cartesian product of complete
transposition graphs \(CT_a\) (\(a\ge 2\)), of order at least \(4\), is
paired \(2\)-disjoint-path-coverable.}

\end{papertheorem}

\hypertarget{input-from-the-disjoint-path-cover-literature}{%
\subsubsection{7.1 Input from the disjoint-path-cover
literature}\label{input-from-the-disjoint-path-cover-literature}}

Throughout this section \(n\) denotes the rank of a transposition-like
graph and \(S,T\) denote terminal sets --- not the paper's global matrix
dimensions or margins.

We use two results of Coleman, Fischberg, Gong, Harrington and Wong
\citep{colemanPairedn1ton1Disjoint2025}. A \textbf{weld} of graphs
\(G_1,\dots,G_w\) of equal order is their disjoint union together with a
perfect matching between every pair \(G_i,G_j\). A
\textbf{transposition-like graph of rank \(1\)} is any
Hamilton-connected or Hamilton-laceable graph; of rank \(r\ge2\), a weld
of \(w\ge r\) transposition-like graphs of rank \(r-1\) of equal order.
Unfolding the recursion produces a welding tower whose bottom layer
consists of rank-\(1\) graphs; these are its \textbf{rank-\(1\) leaves}.
Each \(CT_a\) is a rank-\(a\) example: its \(a\) cosets are copies of
\(CT_{a-1}\), matched by the transpositions moving the last symbol.

\begin{papertheorem}

\textbf{Theorem 7.2} (\citep[Thm.
1.5]{colemanPairedn1ton1Disjoint2025}). \emph{Let \(n \ge 2\) be an
integer and let \(H\) be a bipartite transposition-like graph of rank
\(n\) in which every rank-\(1\) leaf of the welding is a single vertex
or has even order. Then for every choice of \((n-1)\)-subsets
\(S\subseteq V_1\), \(T\subseteq V_2\), \(H\) admits a paired
\((n-1)\)-disjoint path cover.}

\end{papertheorem}

\begin{papertheorem}

\textbf{Theorem 7.3} (\citep[Prop.
1.1(c)]{colemanPairedn1ton1Disjoint2025}). \emph{A balanced bipartite
graph with \(|V_1|=|V_2|\ge k\) (so that \(k\)-subsets of each side
exist) that admits a paired \(k\)-disjoint path cover for every choice
of \(k\)-subsets \(S\subseteq V_1\), \(T\subseteq V_2\) also admits a
paired \(\ell\)-disjoint path cover for every such choice of
\(\ell\)-subsets and every \(1\le\ell\le k\).}

\end{papertheorem}

(This size proviso is implicit in
\citep{colemanPairedn1ton1Disjoint2025}, whose proof extends an
\(\ell\)-demand upward to a \(k\)-demand; without it the \(k\)-subset
hypothesis is vacuously satisfied while the \(\ell\)-conclusion can
fail, so we state it explicitly. It holds wherever the theorem is used
below. The terminal sets \(S\subseteq V_1\) and \(T\subseteq V_2\) below
are \textbf{indexed}, \(S=\{s_1,\ldots,s_k\}\) and
\(T=\{t_1,\ldots,t_k\}\), the \(i\)th path of the cover joining \(s_i\)
to \(t_i\); in that form they are exactly the pairwise-opposite
\(k\)-demands of §2.4, by the within-pair relabeling given there.)

Theorem 7.2 is invoked exactly as stated, and always at rank
\(n \ge 2\). Its only condition on the rank-\(1\) leaves is that each be
a single vertex or an even-order Hamilton-laceable graph. (The size
requirements internal to the rank induction of
\citep{colemanPairedn1ton1Disjoint2025} --- in particular the order
floor \(4n-2\), which appears only in Proposition 1.6 of that paper and
never as a hypothesis of its Theorem 1.5 --- are discharged there by its
own base cases; they are not additional conditions we must check on the
leaves.) In particular Theorem 7.2 applies with leaves that are
themselves products of complete transposition graphs --- the observation
that drives the induction below.

\textbf{The feature we use.} The whole bipartite branch concentrates on
a single feature of Theorem 7.2: that a rank-\(1\) leaf may be
\emph{any} even-order Hamilton-laceable graph, with no further structure
required. We use nothing else about the leaves, and in §7.3 the leaves
are themselves products of complete transposition graphs, whose
even-order and Hamilton-laceability hypotheses §7.3 checks explicitly.

\hypertarget{a-single-block}{%
\subsubsection{7.2 A single block}\label{a-single-block}}

For \(a\ge3\), \(CT_a\) is a rank-\(a\) bipartite transposition-like
graph with single-vertex leaves, so by Theorem 7.2 it admits a paired
\((a-1)\)-disjoint path cover for all pairwise-opposite demands; by
Theorem 7.3 it admits a paired \(2\)-disjoint path cover, i.e.~\(CT_a\)
is paired \(2\)-DPC. (For \(a=2\), \(CT_2=K_2\) has order \(<4\) and is
excluded; co-permutation blocks are isomorphic to \(CT_a\).)

\hypertarget{arbitrary-products}{%
\subsubsection{7.3 Arbitrary products}\label{arbitrary-products}}

\emph{Proof of Theorem 7.1, by induction on the number \(k\) of
factors.}

For \(k=1\) this is §7.2. If all factors are \(CT_2\), the product is a
hypercube \(Q_k\), which is paired \(2\)-DPC for pairwise-opposite
demands by Jo, Park and Chwa \citep[Lem. 1]{joPaired2disjointPath2013}
(a lemma they record from earlier work, valid in every dimension
\(\ge 2\); the one-dimensional \(Q_1=K_2\) is immediate).

Otherwise some factor has rank \(\ge3\); since Cartesian products
commute up to isomorphism, reorder the factors so that one such factor
leads, and write \(G=CT_a\,\square\,B\) with \(a\ge3\) and \(B\) the
product of the remaining \(k-1\) factors. We need only that \(B\) is
even-order and Hamilton-laceable (the rank-\(1\) leaf condition of
Theorem 7.2). If \(|V(B)|\ge 4\) this is the induction hypothesis (\(B\)
is paired \(2\)-DPC, hence Hamilton-laceable by Theorem 7.3 with
\(\ell=1\)); the sole remaining case is \(B=CT_2=K_2\) (one remaining
factor, of rank \(2\)), which is even-order and Hamilton-laceable
directly. Either way \(B\) is balanced bipartite of even order
(Proposition 2.2).

Now \(CT_a\,\square\,B\) is a rank-\(a\) bipartite transposition-like
graph whose rank-\(1\) leaf is \(B\): it is the weld of \(a\) copies of
\(CT_{a-1}\,\square\,B\) (the coset matchings of \(CT_a\), with the
identity on \(B\)), the welding tower bottoming out at
\(CT_1\,\square\,B=B\) (with \(CT_1=K_1\), the one-vertex graph). Since
\(B\) is even-order and Hamilton-laceable, Theorem 7.2 applies and
yields a paired \((a-1)\)-disjoint path cover of \(CT_a\,\square\,B\);
by Theorem 7.3 it admits a paired \(2\)-disjoint path cover.

The two inductions act on different parameters and do not interact: the
inner one (Theorem 7.2) fixes \(B\) as a leaf and recurses on rank,
while the outer one recurses on the number of factors and supplies the
single fact \(B\) must satisfy --- its Hamilton-laceability. \(\square\)

By §7.3 and Proposition 2.2, every balanced bipartite interchange graph
is paired \(2\)-disjoint-path-coverable, which is the Bipartite 2-DPC
Lemma. With the reduction of Sections 3--6, this proves Theorem 1.1.
\(\blacksquare\)

\begin{center}\rule{0.5\linewidth}{0.5pt}\end{center}

\hypertarget{concluding-remarks}{%
\subsection{8. Concluding remarks}\label{concluding-remarks}}

The proof situates Brualdi's problem inside the disjoint-path-cover
theory developed for interconnection networks: the structure theory of
\citep{brualdiCombinatorialMatrixClasses2006} presents the balanced
bipartite interchange graphs as products of complete transposition
graphs, and the transposition-like-graph theorem of
\citep{colemanPairedn1ton1Disjoint2025} then supplies the paired
disjoint path covers that the reduction lifts to all interchange graphs.
The strengthening from Hamiltonicity to maximal Hamiltonicity is what
makes the reduction possible.

Otherwise the interchange-graph literature has studied \(G(R,S)\) mainly
for diameter and connectivity (Brualdi--Manber
\citep{brualdiPrimeInterchangeGraphs1983}; Brualdi \citep[Ch.
6]{brualdiCombinatorialMatrixClasses2006}). Previous Hamiltonicity
results concern special margins: Zhang--Zhang and H. Zhang obtained
Hamilton cycles for certain special margin vectors, Li and Zhang
\citep{liHamiltonicityTypeInterchange1994a} proved Hamiltonicity when
one margin is all-ones \citep[Thm. 2.3,
pp.~109--110]{liHamiltonicityTypeInterchange1994a} and stated
immediately afterward that their method extends to edge-Hamiltonicity
\citep[p.~110]{liHamiltonicityTypeInterchange1994a}, and Arikati and
Peled \citep{arikatiRealizationGraphDegree1999} settled the
majorization-gap-one case in the realization-graph formulation. The
present work extends that program to arbitrary margins, and from
Hamiltonicity to maximal Hamiltonicity. Doing so requires two
ingredients absent from those special structures: the
disjoint-path-cover strengthening of
\citep{colemanPairedn1ton1Disjoint2025, joPaired2disjointPath2013}, used
to thread products of nontrivial blocks, and the pivot construction of
Section 5, used for the indecomposable non-product classes (which
relies, in turn, on the Naddef--Pulleyblank dichotomy: a non-bipartite
matroid base graph is Hamilton-connected \citep{NP1981}). The recent
inputs are the transposition-like-graph theorem
\citep{colemanPairedn1ton1Disjoint2025} and the hypercube-like covers
\citep{joPaired2disjointPath2013}; as the newest and most specialized
results the argument invokes, and the enabling input for the bipartite
branch, they are among the external theorems the formalization proves
from first principles rather than assuming; what the machine-checked
result still takes on faith is the classical combinatorial-matrix and
polyhedral theory of Table 1, together with Alspach's theorem that
Johnson graphs are maximally Hamiltonian. The prime-block classification
itself is due to Brualdi's structure theory. What is new here is the use
of that classification to bring complete-transposition products under
the disjoint-path-cover theorem and thereby close the full bipartite
branch.

A natural question is whether the argument extends to the realization
graphs of arbitrary degree sequences, Problem P59 of
\citep{mutzeCombinatorialGrayCodes2023}, of which Brualdi's problem is
the split-sequence case by the correspondence of Arikati and Peled
\citep{arikatiRealizationGraphDegree1999}. In that generality Hladík and
Fink \citep{hladikRealizationGraphEvery2026} announce a Hamilton path
starting from any prescribed vertex; what maximal Hamiltonicity asks
beyond this is control of \emph{both} endpoints, and it is that
strengthening the present machinery would have to supply. It does not
extend directly, and the obstruction is structural rather than
technical. The engine of Section 5 is that the chosen pivot quotient is
a non-bipartite matroid base graph (Lemmas 5.3, 5.4, and 5.10), so the
Naddef--Pulleyblank theorem applies to it; this rests on the
independence of row and column constraints in the margin setting. For a
general degree sequence the analogous quotient, the family of admissible
neighborhoods of a fixed vertex, fails matroid basis exchange, with
counterexamples on as few as thirteen vertices, so the classical
Hamiltonicity theory of matroid base graphs says nothing about these
quotients, and new structure theory is required. A companion paper in
preparation develops the parts of the present machinery that do
transfer, including the product composition and a strengthening of the
known triangle-free case. Whether those methods deliver the
maximally-Hamiltonian strengthening there remains open.

A second question is whether the conclusion strengthens from
Hamiltonicity to pancyclicity. The dichotomy is forced here as well: a
bipartite graph has no odd cycles, so the statement to aim at is that
\(G(R,S)\) is \textbf{pancyclic} when it is not bipartite, containing a
cycle of every length from \(3\) to \(|V(G)|\), and \textbf{bipancyclic}
when it is, containing a cycle of every even length from \(4\) to
\(|V(G)|\). Computation on small classes is consistent with this.

The quotient is not the obstacle. The present argument asks the pivot
quotient only for Hamilton-connectivity, which the Naddef--Pulleyblank
theorem supplies; a pancyclic conclusion asks it for cycles of
prescribed length through a prescribed edge, and for matroid base graphs
that is known, in a dichotomy of the same shape. Bondy and Ingleton
\citep{BI1976} proved that a matroid basis graph is either a cube, every
edge of which lies on cycles of every even length, or a graph in which
every edge lies on cycles of every length, save that an edge whose
symmetric difference is a component of the matroid misses only the
triangles; Naddef \citep{N1984} extends this to a class of
\((0,1)\)-polyhedra containing the matroid basis polytopes. By Lemmas
5.3 and 5.4 the pivot quotient is exactly such a graph, so on the
quotient the pancyclic input is available where the Hamiltonian one was.

The difficulty lies in the gluing. Section 5 threads the fibers to
produce a single Hamilton path, whereas a pancyclic conclusion requires
the same construction to realize every admissible length, and so to
control how the lengths contributed by the separate fibers combine. That
is the compatibility at the crossings described in the closing note
below, now asked of a whole spectrum rather than of one traversal.
Whether it can be met, or avoided, remains open.

A closing note on method. The reduction developed here dissolves a
global traversal into fiber traversals joined at crossing edges, and
nearly every part of that dissolution --- the block structure, the
counting, the classification of base classes --- yields to systematic
argument; the exception is the requirement that choices made
independently in separate fibers be made \emph{compatibly}, a
synchronization that no local or counting bound supplies and that must
be secured by hand at the crossings. That coupling, rather than any
single lemma, is where the difficulty of the theorem concentrates.

\hypertarget{formal-verification-and-computational-corroboration}{%
\subsection{9. Formal verification and computational
corroboration}\label{formal-verification-and-computational-corroboration}}

Two independent machine checks accompany the proof: a full formalization
in Lean 4, and an exhaustive SAT-based computation on all small classes.
Neither replaces the cited results, which are used as published, and
neither forms part of the human-readable argument of Sections 3--7; they
corroborate it.

\hypertarget{the-formal-development}{%
\subsubsection{The formal development}\label{the-formal-development}}

The released artifact is the companion repository
(\url{https://github.com/jbaggett/brualdi-interchange-lean}) at the
immutable tag \texttt{arxiv-v3}. Its public entry point,
\texttt{brualdi\_lean/TRUST\_SURFACE.md}, gives the final theorem
declarations, the complete axiom trace, a prose-to-Lean theorem map, the
exact statements and source locators for the seven assumed theorems, and
the trust boundary. The tagged tree includes \texttt{lean-toolchain},
\texttt{lakefile.toml}, and \texttt{lake-manifest.json}. From the
repository root the build gate is

\begin{verbatim}
cd brualdi_lean
lake build
test -f .lake/build/lib/lean/BrualdiLean/Sec5.olean
lake env lean CheckAxioms.lean
\end{verbatim}

where \texttt{CheckAxioms.lean} is a two-line file the reader creates in
\texttt{brualdi\_lean/} --- \texttt{import\ BrualdiLean.Ledger} followed
by \texttt{\#print\ axioms\ Brualdi.Ledger.brualdi\_MH} --- and the
expected trace is printed verbatim in the public trust surface. The
internal labels \textbf{CORE} (the Bipartite 2-DPC implication
represented by Lemma 1.2) and \textbf{Brualdi-MH} (the final
maximal-Hamiltonicity theorem) are formal-project labels, not additional
mathematical hypotheses. The composed statement, that every interchange
graph is maximally Hamiltonian, is machine-checked from Lean's
foundations together with exactly seven cited results, and the axiom
trace is mechanically reproducible.

\hypertarget{the-seven-assumed-results}{%
\subsubsection{The seven assumed
results}\label{the-seven-assumed-results}}

The Lean development assumes exactly seven results of the literature as
axioms, listed in Table 1. Two inputs the paper also cites are, by
contrast, proved from first principles within the formalization rather
than assumed: the disjoint-path-cover results of Section 7 and the
Gale--Ryser existence theorem underlying Lemma 5.3. Six of the seven
axioms are named theorems and one is a classical transportation-polytope
dimension count, with the caveat that the complete-transposition product
decomposition (row 3) is a synthesis of Brualdi's Theorem 6.3.4 and the
product structure of his §6.3 rather than a single displayed theorem.
Row 4, pivot-freedom, is our Lemma 5.2.

\begin{longtable}[]{@{}
  >{\raggedright\arraybackslash}p{(\columnwidth - 4\tabcolsep) * \real{0.5769}}
  >{\raggedright\arraybackslash}p{(\columnwidth - 4\tabcolsep) * \real{0.3077}}
  >{\raggedright\arraybackslash}p{(\columnwidth - 4\tabcolsep) * \real{0.1154}}@{}}
\caption{The seven cited results the formalization assumes as axioms;
every other step in the proof is machine-checked from Lean's
foundations. The disjoint-path-cover imports of Section 7 --- Theorem
7.2 with its Proposition 1.6, Theorem 7.3, and the paired hypercube
covers --- are instead proved within the formalization from first
principles, so those citations record provenance rather than assumption,
as Section 7 notes where they occur.}\tabularnewline
\toprule\noalign{}
\begin{minipage}[b]{\linewidth}\raggedright
Cited result
\end{minipage} & \begin{minipage}[b]{\linewidth}\raggedright
Reference
\end{minipage} & \begin{minipage}[b]{\linewidth}\raggedright
Used in
\end{minipage} \\
\midrule\noalign{}
\endfirsthead
\toprule\noalign{}
\begin{minipage}[b]{\linewidth}\raggedright
Cited result
\end{minipage} & \begin{minipage}[b]{\linewidth}\raggedright
Reference
\end{minipage} & \begin{minipage}[b]{\linewidth}\raggedright
Used in
\end{minipage} \\
\midrule\noalign{}
\endhead
\bottomrule\noalign{}
\endlastfoot
The interchange graph \(G(R,S)\) is connected & Ryser
\citep{ryserCombinatorialPropertiesMatrices1957}; Brualdi--Manber
\citep{brualdiPrimeInterchangeGraphs1983} & §2 \\
A non-bipartite interchange graph contains a triangle & Brualdi
\citep[Thm. 6.3.4]{brualdiCombinatorialMatrixClasses2006} & §2 \\
A balanced bipartite interchange graph is a Cartesian product of
complete-transposition graphs & Brualdi
\citep[§6.3]{brualdiCombinatorialMatrixClasses2006} & §2 \\
In a nonempty active prime class every row--column cell is variable
(pivot-freedom) & Brualdi--Manber \citep[Thm.
9]{brualdiPrimeInterchangeGraphs1983} & §5 \\
A non-bipartite matroid base-polytope graph is Hamilton-connected &
Naddef--Pulleyblank \citep[Thm. 3.3.1 and Cor. 3.3.2]{NP1981} & §5 \\
Johnson graphs are maximally Hamiltonian & Alspach \citep[Thm.
1.1]{alspachJohnsonGraphsAre2012} & §6 \\
A transportation-polytope dimension count & Brualdi--Hartfiel--Hwang
\citep[§9.13]{brualdiCombinatorialMatrixClasses2006} & §6 \\
\end{longtable}

\hypertarget{correspondence-with-the-printed-proof}{%
\subsubsection{Correspondence with the printed
proof}\label{correspondence-with-the-printed-proof}}

The companion's proofs are the arguments printed here; an adversarial
correspondence audit confirmed the match claim by claim. Proposition 4.1
(the walk device, its Claim, and the boundary pass), Lemma 5.3 (on top
of a self-contained mechanization of the Gale--Ryser existence theorem),
Lemma 5.7, and Lemma 5.10 (including the shiftedness of its pattern
matroid) are machine-checked as printed, each of Proposition 4.1, Lemma
5.3, and Lemma 5.10 also carrying an independent machine-checked
alternate proof --- for Proposition 4.1 an absorber construction of the
same covers; for Lemma 5.3 the minimal-pair argument; for Lemma 5.10 two
alternates, a compression by rank-sum minimization and a shiftedness
proof by residual unit transfer --- and Lemma 5.7 an independent
machine-checked proof of the avoid form the gluing consumes, which
handles all \(D\le 1\) without the two exclusion steps. Theorem 7.1 is
machine-checked in a form that subsumes its statement (arbitrary factor
orders, with no order hypothesis: below four vertices the property is
vacuous). This separates the kernel-checked deduction from its
literature assumptions and from the external SAT tests.

\hypertarget{computational-corroboration}{%
\subsubsection{Computational
corroboration}\label{computational-corroboration}}

Independently of the proof, we verified Theorem 1.1 directly by computer
for all small margins. For every canonical class \((R,S)\) up to 7×7
small enough to enumerate (144,389 isomorphism classes, complete through
\(|V|=35\)), the interchange graph \(G(R,S)\) was confirmed maximally
Hamiltonian by deciding each required Hamilton-path instance (2,555,559
in total) with a SAT solver, with no exception. This finite check
corroborates the result independently of the argument: the SAT sweep
tests its instances, the mechanization its deductions.

\hypertarget{acknowledgments}{%
\subsection{Acknowledgments}\label{acknowledgments}}

The computational search and verification programs and the Lean 4
formalization (Section 9) were developed with AI-assisted tools
(Anthropic Claude; OpenAI GPT/Codex) under author direction. All
mathematical content was verified independently of any model's
assertions: the argument is machine-checked in Lean 4 (Section 9) and
validated computationally on all small classes. The authors take full
responsibility for the content. This use of AI-assisted tools follows
the principles of the Leiden Declaration on Artificial Intelligence and
Mathematics \citep{leidenDeclarationAI2026}, including its call for the
formal verification of AI-assisted results where appropriate --- met
here by the Lean 4 machine-checking --- and for correctness and
responsibility to remain with the human authors.

\hypertarget{funding}{%
\subsection{Funding}\label{funding}}

This research did not receive any specific grant from funding agencies
in the public, commercial, or not-for-profit sectors.

\hypertarget{data-availability}{%
\subsection{Data availability}\label{data-availability}}

No datasets were generated or analyzed for this article. The Lean 4
formalization and verification artifacts supporting its findings are
openly available in the public companion repository at
\url{https://github.com/jbaggett/brualdi-interchange-lean}.

\hypertarget{declaration-of-generative-ai-and-ai-assisted-technologies-in-the-writing-process}{%
\subsection{Declaration of generative AI and AI-assisted technologies in
the writing
process}\label{declaration-of-generative-ai-and-ai-assisted-technologies-in-the-writing-process}}

During the preparation of this work the authors used Anthropic Claude
and OpenAI GPT/Codex in order to draft and edit prose. After using these
tools, the authors reviewed and edited the content as needed and take
full responsibility for the content of the publication. The tools' role
in developing the mathematics itself --- proof search, the computational
programs, and the Lean 4 formalization --- is described in the
Acknowledgments.

\hypertarget{references-section}{%
\subsection{References}\label{references-section}}

  \bibliography{references.bib}

\end{document}